%% file: tropicaljinvariant.tex
\documentclass [10pt,reqno]{amsart}

\usepackage {hannah}
\usepackage{pstricks}
\usepackage{pst-plot}
\usepackage {color,graphicx}
\usepackage {multirow}
\usepackage{texdraw}
\usepackage [latin1]{inputenc}
\usepackage[all]{xy}

\def\ol#1{\overline{#1}}
\DeclareMathOperator {\Int}{Int}
\DeclareMathOperator {\weight}{weight}
\DeclareMathOperator {\cyclelength}{cl}
\DeclareMathOperator {\val}{val}
\DeclareMathOperator {\lc}{lc}
\DeclareMathOperator {\tini}{t-in}
\DeclareMathOperator {\LC}{LC}
\DeclareMathOperator {\trop}{trop}
\DeclareMathOperator {\Trop}{Trop}
\DeclareMathOperator {\UH}{UH}
\DeclareMathOperator {\PGl}{PGl}
\DeclareMathOperator {\Conv}{Conv}
\DeclareMathOperator {\vol}{vol}

\newcommand{\DET}[2]{\det\left(\begin{array}{#1} #2\end{array}\right)}
\newcommand{\ca}{{\mathcal{A}}}

\newcommand{\cx}{{\mathcal{X}}}
\newcommand{\cu}{{\mathcal{U}}}
\newcommand{\cv}{{\mathcal{V}}}
\newcommand{\cf}{{\mathcal{F}}}

\CompileMatrices     

\title [The $j$-invariant of a plane tropical cubic]{The $j$-invariant of a plane tropical cubic}

\author {Eric Katz, Hannah Markwig, Thomas Markwig}
\address {Eric Katz, Department of Mathematics, The University of Texas
  at Austin, 1 University Station, C1200, Austin, TX 78712}
\email {eekatz@math.duke.edu}
\address{Hannah Markwig, University of Michigan, Department of Mathematics, 2074 East Hall, 530 Church Street, Ann Arbor, MI 48109-1043}
\email{markwig@umich.edu}
\address {Thomas Markwig, Fachbereich Mathematik, Technische Universit\"at Kaiserslautern, Postfach 
3049, 67653 Kaiserslautern,
Germany}
\email {keilen@mathematik.uni-kl.de}
\thanks{The third author would like to thank the Institute for Mathematics and its Applications in 
Minneapolis for hospitality}

\begin {document}

\begin {abstract}  
  Several results in tropical geometry have related the $j$-invariant
  of an algebraic plane curve of genus one to the cycle length of a
  tropical curve of genus one. In this paper, we 
  prove that for a plane cubic over the field of Puiseux
  series the negative of the generic valuation of the $j$-invariant is
  equal to the cycle length of the tropicalization of the curve, if
  there is a cycle at all.  
\end {abstract}

\maketitle

\section{Introduction}

Tropical geometry is a new and active field of research. Roughly, its
main idea is to replace algebraic varieties by piece-wise linear
objects called tropical varieties. These objects may be easier to deal
with and new methods from combinatorics can be used to handle
them. A lot of work has been done to ``translate'' terms and
definitions to the tropical world. Sometimes a translation is
justified by the appropriate use of the new tropical object rather
than by an argument why this is the correct tropicalization. This is
the case for the $j$-invariant of an elliptic curve. Many results
predict that the ``tropical $j$-invariant'' of a tropical curve of
genus one is its cycle length.

The $j$-invariant is an invariant which coincides for two smooth
elliptic curves over an algebraically closed fields if and only if
they are isomorphic. In \cite{Mi06}, 
isomorphisms (``equivalences'') between abstract tropical curves are
defined, and two elliptic abstract tropical curves are equivalent if and
only if they have the same cycle length. Thus the cycle length plays
the same role in the tropical setting as the $j$-invariant does in the
algebraic setting. 
Also, the possibility to define a group law on the cycle of a tropical
curve (see \cite{Vig04}) using distances indicates the importance of
the cycle and its length. 
Furthermore, the numbers of tropical curves with fixed cycle length
are in correspondence to the numbers of curves with fixed
$j$-invariant (see \cite{KM06}). But there is also a result
which suggests that the cycle length actually might be the correct
tropicalization  of the $j$-invariant (as introduced on p.\ \pageref{p-trop}).
This is a byproduct of the proof of \cite[Theorem 5.4.1]{Spe05}. It
says that given a tropical curve $\mathcal{C}$ with cycle length $l$
and a Puiseux series $j$ of valuation $-l$, we can embed the elliptic
curve with $j$-invariant $j$ such that its tropicalization is equal to
$\mathcal{C}$.  

All these results indicate that the ``tropical $j$-invariant'' should
be defined as the cycle length (respectively, its negative). The aim
of this paper is to  
show that for a plane cubic  the
$j$-invariant really tropicalizes to the negative of the cycle
length. 

More precisely, we define plane cubic curves over the field of Puiseux
series $\K=\C\{\{t\}\}$ and use the 
valuation map to tropicalize them (see p.\ \pageref{p-trop}).  
The $j$-invariant of an elliptic curve over the Puiseux series is a
Puiseux series itself. 
Our main theorem (see Theorem \ref{thm-main}) is that the negative of
the cycle length
of the tropicalization of a smooth cubic curve in $\P_{\K}^2$
(assuming it has a cycle -- see Definition \ref{def-cycle})
is always equal to the generic valuation of the
$j$-invariant (see Definition \ref{def-genericvaluation}) and it is
actually equal to 
the valuation of the $j$-invariant itself  if no
terms in the $j$-invariant cancel -- which generically is the case. 
A corollary (see Corollary \ref{cor-nocycle}) of this theorem is that if an elliptic
curve has a $j$-invariant with a positive valuation, then its
tropicalization does not have a cycle. 

There is an intriguing similarity to bad reduction of elliptic curves over
discrete valuation rings. Firstly, only elliptic curves whose
$j$-invariant has a negative valuation have bad reduction, and
secondly the negative of this valuation is then the ``cycle length'' of the
special fiber in the N\'eron model of the curve in the sense that it
is the number of projective lines forming the cycle.
We are working on better understanding of the
connection of these two results. Moreover, in a forthcoming paper we
will show that the main result of this paper can be generalized to
smooth elliptic curves on toric surfaces other than the projective
plane. 

This paper is organized as follows. In Section \ref{sec-jinv} we
recall the definition of the $j$-invariant of a plane cubic as a
rational function in the cubic's coefficients. Its
denominator is the discriminant of the cubic. Moreover, we observe
that the generic valuation (see Definition \ref{def-genericvaluation}) 
of the $j$-invariant is a piece-wise linear function. In
Section \ref{sec-trop} we recall basic definitions of tropical
geometry and show that the function ``cycle length'' is piece-wise
linear as well. 
The main theorem is stated in Section \ref{sec-thm}. As we know
already that the two functions we want to compare are piece-wise
linear the proof consists of two main steps: first we compare certain
domains of linearity, then we compare the two linear functions on each
domain. We present two proofs since we believe that each of them is
interesting on its own. For the first proof, we choose as domains of
linearity the union of $\Delta$-equivalent cones of the secondary fan
of $\ca_3=\{(i,j)\;|\;0\leq i,j,i+j\leq 3\}$ (i.e.\ cones of the
Gr\"obner fan of the discriminant). The 
comparison of the two linear functions ``generic valuation of the
$j$-invariant'' and ``cycle length'' on each such cone
is done in Section \ref{sec-delta}. In the second proof, we
choose smaller domains of linearity --- cones of the secondary fan of
$\ca_3$. To 
compare the two linear functions, we have to classify the rays of the
secondary fan of $\ca_3$ and compare them on each ray. Section \ref{sec-rays} is
concerned with this classification of the rays of the secondary fan of
$\ca_3$. In Section 
\ref{sec-discriminant} we give an alternative proof of the fact that
for an arbitrary convex lattice polytope the Gr\"obner fan of the
discriminant is a
coarsening of the secondary fan, and
in Section \ref{sec-num} we study the numerator of the
$j$-invariant. These 
two sections are important to understand the domains of linearity of
the function ``generic valuation of the $j$-invariant''.  

Parts of our proofs and many examples rely on computations performed
using \texttt{polymake} \cite{Pol97}, TOPCOM \cite{topcom} and \textsc{Singular}
\cite{GPS05}. 
The \textsc{Singular} code that we used for this is contained in the
library \texttt{jinvariant.lib} (see \cite{KMM07a}) and it is  
available via the URL 
\begin{center}
  http://www.mathematik.uni-kl.de/\textasciitilde keilen/en/jinvariant.html.
\end{center}
More detailed
explanations on how to use the code can be found there. 

The authors would like to acknowledge Vladimir Berkovich, Jordan
Ellenberg, Bjorn Poonen, David Speyer, Charles Staats, Bernd Sturmfels
and John Voigt for valuable discussions. 


\section{The j-invariant and its valuation}\label{sec-jinv}

Since every smooth elliptic curve can be embedded into the projective
plane as a cubic it makes sense to start the investigation of smooth elliptic
curves and their tropicalizations by studying cubics in the plane. And
since the tropicalization of a curve highly depends on its embedding
we want to consider all such planar embeddings of a given curve at the
same time. That is, we start with a non-zero homogeneous polynomial 
\begin{displaymath}
  f=\sum_{i+j=0}^3 a_{ij}x^iy^jz^{3-i-j}
\end{displaymath}
of degree $3$ as input data. Here the
coefficients $a_{ij}$ are thought of as elements of the field 
\begin{displaymath}
  \K=\bigcup_{N=1}^\infty\Quot\Big(\C\big[\big[t^{\frac{1}{N}}\big]\big]\Big)
  =\left\{\sum_{\nu=m}^\infty c_\nu\cdot
    t^\frac{\nu}{N}\;\Big|\;
    c_\nu\in\C, N\in\Z_{>0}, m\in\Z\right\}
\end{displaymath}
of Puiseux series over $\C$, and we consider the algebraic curve
$C=V(f)\subset\P_{\K}^2$ defined by $f$. Since the tropicalization of
$C$ only depends on the points of $C$ in the torus
$(\K^*)^2=\{(x:y:1)\;|\;x\not=0\not=y\}$ (see Section \ref{sec-trop}) we
may as well replace $f$ by its affine equation
\begin{displaymath}
  f=\sum_{i+j=0}^3 a_{ij}x^iy^j,
\end{displaymath}
and we will do so for the remaining part of the paper -- always
keeping in mind that via homogenization $f$ defines a cubic in the
projective plane. 

In our investigation the Newton polytope of $f$, the triangle $Q_3$
with endpoints $(0,0)$, $(0,3)$ and $(3,0)$, plays an important
role. We denote by $\ca_3:=Q_3\cap \Z^2$ its integer points. That way
we can write the equation for $f$ as $f=\sum_{(i,j)\in\ca_3}a_{ij}x^iy^j$.

\begin{figure}[h]
  \begin{texdraw}
    \drawdim cm \relunitscale 0.4 \linewd 0.05 \lpatt (1 0) \setgray 0.8 
    \move (0 0) \lvec (3 0) \lvec (0 3) \lvec (0 0) \lfill f:0.8
    \move (0 0) \fcir f:0 r:0.1 
    \move (0 1) \fcir f:0 r:0.1 
    \move (0 2) \fcir f:0 r:0.1 
    \move (0 3) \fcir f:0 r:0.1 
    \move (1 0) \fcir f:0 r:0.1 
    \move (1 1) \fcir f:0 r:0.1 
    \move (1 2) \fcir f:0 r:0.1
    \move (2 0) \fcir f:0 r:0.1 
    \move (2 1) \fcir f:0 r:0.1 
    \move (3 0) \fcir f:0 r:0.1 
  \end{texdraw}\centering
  \caption{$Q_3$ and $\ca_3$}\label{fig-q3}
\end{figure}

If the curve $C$ is smooth it is fixed up to isomorphism by a single
invariant $j(C)=j(f)\in \K$, the \emph{$j$-invariant}, which can be computed as
a rational function, say
\begin{displaymath}
  j(f)=\frac{A}{\Delta}  
\end{displaymath}
in the coefficients $a_{ij}$ of $f$ -- now thought of as indeterminates
-- where $A,\Delta\in\Q[a_{i,j}\;|\;(i,j)\in\ca_3]$ are
homogeneous polynomials of degree $12$. Since we quite frequently need
to refer to polynomials in the $a_{ij}$ we introduce the convention
$\underline{a}=(a_{ij}\;|\;(i,j)\in\ca_3)$ and if
$\omega\in\N^{\ca_3}$ is a multi index then $\underline{a}^\omega=\prod_{(i,j)\in\ca_3}
a_{ij}^{\omega_{ij}}$. 
The denominator $\Delta$
is actually the \emph{discriminant} of $f$ (see \cite{GKZ}). 

The field $\K$ of Puiseux series comes in a natural way with a
\emph{valuation}, namely
\begin{displaymath}
  \val:\K^*\rightarrow\Q:\sum_{\nu=m}^\infty c_\nu\cdot
  t^\frac{\nu}{N}\mapsto \min\left\{\frac{\nu}{N}\;\Big|\; c_\nu\not=0\right\},
\end{displaymath}
and we may extend the valuation to $\K$ by $\val(0)=\infty$.  
If $k=\sum_{\nu=m}^\infty c_\nu\cdot t^\frac{\nu}{N}\in\K$ with
$\val(k)=\frac{m}{N}$ then we call $\lc(k):=c_m$ the \emph{leading coefficient}
of the formal power series $k$. We
sometimes call $\val(k)$ the \emph{tropicalization} of $k$.\label{p-trop}

Throughout the paper we will treat polynomials in the variables
$(x,y)$ as well as in the variables
$\underline{a}=(a_{ij}\;|\;(i,j)\in\ca_3)$, and many results will be
derived for both cases at the same time. We therefore want to
introduce a unifying notation which we use whenever we work with
polynomials with an arbitrary set of variables.

\begin{notation}\label{not-x}
  We set $\underline{x}=(x_\lambda\;|\;\lambda\in\Lambda)$ where
  $\Lambda$ is some finite index set. If
  $\omega=(\omega_\lambda\;|\;\lambda\in\Lambda)\in\N^\Lambda$ then we
  use the usual multi index 
  \begin{displaymath}
    \underline{x}^\omega=\prod_{\lambda\in\Lambda}x_\lambda^{\omega_\lambda},
  \end{displaymath}
  and if $\ca\subset\N^\Lambda$ is finite and $h_\omega\in\K$ for
  $\omega\in\ca$ then
  \begin{displaymath}
    h=\sum_{\omega\in\ca}h_\omega\cdot\underline{x}^\omega\in\K[\underline{x}]
  \end{displaymath}
  is a polynomial over $\K$.
  Later on we will sometimes need Laurent polynomials instead of
  polynomials and we then allow negative exponents.
\end{notation}

It is our aim to study the valuation of the $j$-invariant of the curve
$C$ given by fixing values for the $a_{ij}$. It would be nice if it
only depended on their valuations, but this is only true if the
leading coefficients of the $a_{ij}$ are sufficiently general. We
therefore introduce the notion of \emph{generic valuation} for a
polynomial like $A$ or $\Delta$. In order to do so we need
$t$-initial forms, a concept which is needed in a broader context
further down.

Given a polynomial $0\not=h=\sum_{\omega\in\ca} h_\omega\cdot
\underline{x}^\omega\in\K[\underline{x}]$ and
a point $v\in\R^{\Lambda}$, we define the
\emph{$v$-weight}
\begin{displaymath}
  \weight_v(h)=\min\{\val(h_\omega)+v\cdot \omega\;|\;h_\omega\not=0\}
\end{displaymath}
of $h$, where $v\cdot
\omega=\sum_{\lambda\in\Lambda}v_\lambda\cdot\omega_\lambda$, and the
\emph{$t$-initial form}
\begin{displaymath}
  \tini_v(h)=\sum_{\val(h_\omega)+v\cdot\omega=\weight_v(h)}
  \lc(h_\omega)\cdot\underline{x}^\omega
\end{displaymath}
of $h$ with respect to $v$.

\begin{definition}\label{def-genericvaluation}
  If $0\not=h=\sum_{\omega\in\ca}h_\omega\cdot\underline{x}^\omega\in
  \Q[\underline{x}]\subset\K[\underline{x}]$ has \emph{constant
    rational coefficients} and $v\in\R^\Lambda$, then
  $\val_v(h):=\weight_v(h)=\min\{v\cdot \omega\;|\;h_\omega\not=0\}$ is
  called the \emph{generic valuation} of $h$ at $v$.
\end{definition}

\begin{remark}
  If $0\not=h\in \Q[\underline{x}]$, then all $v\in \R^\Lambda$
  such that the $t$-initial forms $\tini_v(h)$ coincide with each other  form
  a cone and the collection of 
  these cones forms the \emph{Gr\"obner 
    fan} of $h$ (see e.g.\  \cite[Chap.\ 2]{Stu96}). 
  Top-dimensional cones correspond to $t$-initial forms which are
  monomials.   
\end{remark}

\begin{lemma}\label{lem-genvallinear}
  If $0\not=h\in\Q[\underline{x}]$, then the function 
  \begin{displaymath} \val_\cdot(h):\R^\Lambda\rightarrow \R:v\mapsto \val_v(h)\end{displaymath}
  is piece-wise linear, and it is linear on a top-dimensional cone of the Gr\"obner fan of $h$.
  Moreover, if $v\in \R^\Lambda$ is in the interior of a top-dimensional cone of the Gr\"obner
  fan of $h$, then $\val_v(h)=\val(h(y))$ for any $y\in (\K^*)^\Lambda$
  with $\val(y)=v$.  
\end{lemma}
\begin{proof}
  Obviously $\val_\cdot(h)$ is piece-wise linear.
  If $v$ is in the interior of a top-dimensional cone of the Gr\"obner
  fan of $h$, $\tini_v(h)=\underline{x}^{\tilde{\omega}}$ is a
  monomial. Therefore the minimum $v\cdot \omega$ in the definition of
  $\val_v(h)$ is attained for only one term, namely for
  $h_{\tilde{\omega}}\cdot \underline{x}^{\tilde{\omega}}$, and
  $\val_\cdot(h):v\mapsto \val_v(h)=v\cdot\tilde{\omega}$ is linear.  
  Furthermore, if $\val(y)=v$ then the terms of $h(y)$ of lowest order in $t$ come from
  those terms of $h$ for which $v\cdot \omega$ is minimal. Thus
  $\val(h(y))=v\cdot\tilde{\omega}=\val_v(h)$. 
\end{proof}

Note that if $y$ is not in the interior of a top-dimensional cone of
the Gr\"obner fan of $h$, then the terms of lowest order in $t$ in
$h(y)$ might cancel and we cannot predict the valuation of $h(y)$. The
generic valuation is, however, the valuation of $h(y)$ under the assumption that
the terms of lowest order in $t$ do not cancel. 

\begin{definition} 
  With the above notation we define the \emph{generic valuation of the $j$-invariant at
    $u\in\R^{\ca_3}$} as $\val_u(j):=\val_u(A)-\val_u(\Delta)$.
\end{definition}

\begin{remark}
  From Lemma \ref{lem-genvallinear} it follows that 
  \begin{displaymath}
    \val_{\cdot}(j):\R^{\ca_3}\longrightarrow\R:u\mapsto\val_u(j)
  \end{displaymath}
  is a
  piece-wise linear function which is linear on intersections $D\cap D'$
  of a top-dimensional cone $D$ of the Gr\"obner fan of $A$ and a
  top-dimensional cone $D'$ of the Gr\"obner fan of $\Delta$. For $u$ in
  the open interior of $D\cap D'$, $\val_u(j)=\val\big(j(f)\big)$ for any
  $f=\sum_{(i,j)\in \ca_3}a_{ij}x^iy^j$ with $\val(a_{ij})=u_{ij}$. 
\end{remark}


\section{Tropicalizations and the cycle length of a plane tropical cubic}\label{sec-trop}

In this section we will study the tropicalization of plane cubics as
well as the tropicalization of the varieties defined by $A$ or
$\Delta$ in $(\K^*)^{\ca_3}$. We therefore start again using the 
general notation \ref{not-x}.

\begin{definition}
  If $h\in\K[\underline{x}]$ then the
  \emph{tropicalization of} $V(h)=\{p\in\K^\Lambda\;|\;h(p)=0\}$,
  \begin{displaymath}
    \Trop\big(V(h)\big)=\overline{\val\big(V(h)\cap
      (\K^*)^\Lambda)}\subseteq\R^\Lambda,
  \end{displaymath}
  is
  the closure of $\val\big(V(h)\cap (\K^*)^\Lambda\big)$ with respect to
  the Euclidean topology in $\R^\Lambda$, where by abuse of notation
  \begin{displaymath}
    \val:(\K^*)^\Lambda\longrightarrow\Q^\Lambda:(k_\lambda\;|\;\lambda\in\Lambda)\mapsto 
    \big(\val(k_\lambda)\;\big|\;\lambda\in\Lambda\big)
  \end{displaymath}
  denotes the Cartesian product of the valuation map from Section \ref{sec-jinv}.
\end{definition}

This definition is not too helpful when it comes down to actually computing
tropical varieties. There it is better to consider tropical
polynomials.

\begin{definition}
  For a \emph{tropical polynomial} $F=\min\{u_\omega+v\cdot
  \omega\;|\;\omega\in\ca\}$, with $\ca\subset\N^\Lambda$ finite and
  $u_\omega\in\R$, (see 
  e.g.\ \cite{RST03}) we define the \emph{tropical hypersurface} associated to
  $F$ to be the locus of non-differentiability of the function
  $$\R^\Lambda\longrightarrow \R:v\mapsto
  \min_{\omega\in\ca}\{u_\omega+v\cdot\omega\}$$ (i.e.\ the
  locus where the minimum is attained by at least two terms). We call
  the convex hull of $\ca$ the
  \emph{Newton polytope} of $F$ respectively of the tropical
  hypersurface defined by $F$. If $\#\Lambda=2$ then we call a
  tropical hypersurface simply a \emph{tropical curve}.
\end{definition}

\begin{remark}
  Given
  $h=\sum_{\omega\in\ca}h_\omega\cdot\underline{x}^\omega\in\K[\underline{x}]$
  we define the \emph{tropicalization of} $h$ to be the tropical
  polynomial
  \begin{displaymath}
    \trop(h)=\min\big\{\val(h_\omega)+v\cdot \omega\;\big|\;\omega\in\ca\big\}.  
  \end{displaymath}
  Then 
  by Kapranov's Theorem (see \cite[Theorem~2.1.1]{EKL04}), $\Trop\big(V(h)\big)$ is
  equal to the tropical hypersurface associated to $\trop(h)$.

  Moreover, this is obviously the collection of $v$ for which
  $\tini_v(h)$ is not a monomial. In particular, if $h\in\Q[\underline{x}]$ has constant
  rational coefficients then it is the codimension-one skeleton of the
  Gr\"obner fan.

  We will use this definition mainly in two different settings:
  \begin{enumerate}
  \item For a plane cubic $f=\sum_{(i,j)\in\ca_3}a_{ij}x^iy^j$ we will
    consider the corresponding \emph{plane tropical cubic}
    $\Trop\big(V(f)\big)\subset\R^2$ given by the 
    tropical polynomial 
    $\min_{(i,j)\in\ca_3}\{u_{ij}+ix+jy\}$ where $\val(a_{ij})=u_{ij}$.
  \item For the numerator and denominator of the $j$-invariant $A$ and
    $\Delta$ which are polynomials in $\Q[\underline{a}]$ we will consider
    their tropicalizations $\Trop\big(V(A)\big)$ and
    $\Trop\big(V(\Delta)\big)$ in $\R^{\ca_3}$. The latter 
    was recently studied in \cite{DFS05}.  As $A$ and $\Delta$ have
    constant rational coefficients these two tropical hypersurfaces
    are equal to the codimension-one skeletons of the Gr\"obner fans of $A$
    respectively $\Delta$. 
  \end{enumerate}
\end{remark}

Tropical hypersurfaces as defined above are dual to certain marked
subdivisions of $\ca$. Let us recall the necessary facts from
\cite{GKZ}, still using Notation \ref{not-x}. 

\begin{definition} 
  A {\em marked polytope} is a pair $(Q,\ca)$ where
  $Q\subset\R^\Lambda$ is a convex lattice polytope and $\ca\subseteq Q\cap\Z^\Lambda$
  contains at least the vertices of $Q$. The set $\ca$ is said to
  be the set of {\em marked lattice points}. 
\end{definition}

The Newton polytope $(Q_3,\ca_3)$ as shown in Figure \ref{fig-q3} is
a marked polytope.

\begin{definition} 
  Let $(Q,\ca)$ be a marked polytope in $\R^\Lambda$ with $\dim(Q)=\#\Lambda$.
  A \emph{marked subdivision} of $(Q,\ca)$
  is a finite family of marked polytopes $\{(Q_i,\ca_i)\;|\;i=1,\ldots,k\}$ such that
  \begin{enumerate}
  \item $(Q_i,\ca_i)$ is a marked polytope with
    $\dim(Q_i)=\#\Lambda$ for $i=1,\ldots,k$,
  \item $Q=\bigcup_{i=1}^k Q_i$ is a subdivision of $Q$, i.e.\
    $Q_i\cap Q_j$ is a face (possibly empty) of $Q_i$ and of $Q_j$
    for all $i,j=1,\ldots,k$, 
  \item $\ca_i\subset \ca$ for $i=1,\ldots,k$, and
  \item $\ca_i\cap(Q_i\cap Q_j)=\ca_j\cap (Q_i\cap Q_j)$ for all $i,j=1,\ldots,k$.
  \end{enumerate}
\end{definition}

We do not mandate that $\bigcup_{i=1}^k \ca_i=\ca$. Example
\ref{ex-markedsubdivision} shows an example of a marked subdivision
of $(Q_3,\ca_3)$.

\begin{definition} 
  Let $S=\{(Q_i,A_i)\;|\;i=1,\ldots,k\}$ and
  $S'=\{(Q'_j,A'_j)\;j=1,\ldots,k'\}$ be subdivisions of 
  $(Q,\ca)$.  We say that $S$ \emph{refines} $S'$ is for all $j=1,\ldots,k'$, the
  collection of $(Q_i,\ca_i)$ so that $Q_i\subseteq Q'_j$ is a
  marked subdivision of $(Q'_j,\ca'_j)$.
\end{definition}

Figure \ref{fig-gcl} shows two marked subdivisions of $(Q_3,\ca_3)$
where the right one is a refinement of the left one.

\begin{remark}
  Let us note that the \emph{coarsest subdivision} of $(Q,\ca)$ is $\{(Q,\ca)\}$.
  
  More interesting examples of subdivisions are the so called
  \emph{regular} or \emph{coherent subdivisions}.
  Given $\psi\in\R^\ca$, we can associate a subdivision as follows.    
  Let the \emph{upper hull} $\UH(\psi)$ of $\psi$ be the convex hull
  of the subset of $\R^\Lambda\times\R$ given 
  by
  \[S=\{(\omega,a)|\omega\in\ca,a\geq\psi(\omega)\}.\]  
  The lower faces
  of $\UH(\psi)$ project to $Q$ giving a subdivision of $Q$.
  Define the \emph{lower convexity} $\LC(\psi)$ of $\psi$ to be the
  function $\LC(\psi):Q\rightarrow\R$ whose graph is the lower
  faces of $\UH(\psi)$.  This is a convex function.  The faces of the
  subdivision are the maximal domains of linearity of $\LC(\psi)$.  For
  a face $Q_i$, define the marked set  
  \[\ca_i=\{\omega\in Q_i\cap\ca \;|\; \psi(\omega)=\LC(\psi)(\omega)\},\]
  the set of all points of $\ca$ in $Q_i$ that lie on the lower faces
  of the upper hull.  Marked subdivisions that arise in this fashion
  are said to be {\em regular}.  

  The set of all $\psi\in\R^{\ca}$ that
  induce the same regular marked subdivision form an open polyhedral cone in
  $\R^\ca$.  All such cones form the \emph{secondary fan} of $\ca$.  The secondary fan
  is the normal fan to a polytope, the \emph{secondary polytope}. 
  The poset of cones in the secondary fan is isomorphic to the poset
  of regular marked subdivisions under refinement.  \emph{Top dimensional
    cones of the secondary fan} of $\ca$ therefore correspond to the
  finest subdivisions where each 
  $Q_i$ is a  $\#\Lambda$-dimensional simplex and $\ca_i$ is the set of
  vertices of $Q_i$. (See \cite[Chap.\ 7]{GKZ}.)
 
  \begin{remark}\label{linearityspace} Note that the minimal cone of the secondary fan is a
  $\#\Lambda+1$-dimensional space $L$ and is given by all $\psi\in\R^\ca$ of
  the form
  \[\psi:\ca\rightarrow\R:\omega\mapsto
  a+v\cdot\omega\]
  with $a\in\R$ and $v\in\R^\Lambda$. 
  All these functions induce the \emph{coarsest subdivision}.
  We call $L$ the
  \emph{linearity space} of the secondary fan of $\ca_3$.
  Every cone of the secondary fan of $\ca_3$ contains its linearity space.
  Therefore, we may consider the secondary fan as living
  in $\R^\ca/L$.  When we speak of \emph{rays} in the secondary fan, we mean
  $n+2$-dimensional cones in the fan containing $L$. 
  \end{remark}

  The secondary fan of $\ca_3$ is an important object because we
  will see that the ``cycle length'' as a function is linear on each
  top-dimensional cone of the secondary fan. Since we have already
  seen that the valuation of the 
  $j$-invariant is linear on each cone of the common refinement of
  the Gr\"obner fans of $A$ and $\Delta$ our strategy will be to
  compare the secondary fan with these two Gr\"obner fans. 
\end{remark}

\begin{remark}\label{rem-direction}
  The \emph{duality} of tropical hypersurfaces in $\R^\Lambda$ with Newton
  polytope $Q_H\subset\R^\Lambda$ and regular marked subdivisions of
  $(Q_H,\ca_H)$ with $\ca_H=Q_H\cap\Z^\Lambda$ is set up as follows.
  Given a tropical polynomial
  $H=\min_{\omega\in\ca_H}\{u_\omega+v\cdot\omega\}$ we associate to $H$
  the regular marked subdivision induced by 
  \begin{displaymath}
    \psi_H:\ca_H\longrightarrow\R:\omega\mapsto u_\omega.
  \end{displaymath}
  This marked subdivision is then dual to the tropical hypersurface
  defined by $H$ in the sense of \cite[Prop.\ 3.11]{Mi03}, which in
  its full generality is rather technical.
  However, for the cases we are interested in it can easily be
  described. 

  If $H$ defines a plane tropical curve $\mathcal{C}_H$ in $\R^2$ then each marked
  polytope of the subdivision of $(Q_H,\ca_H)$ is dual to a
  vertex of $\mathcal{C}_H$, and each facet of a marked polytope is dual
  to an edge of $\mathcal{C}_H$. Moreover, if the facet, say $F$, has end points
  $(x_1,y_1)$ and $(x_2,y_2)$ then the \emph{direction vector} $v(E)$ of the
  dual edge $E$ in $\mathcal{C}_H$ is defined (up to sign) as
  \begin{displaymath}
    v(E)=(y_2-y_1,x_1-x_2)
  \end{displaymath}
  and points in the direction of $E$.
  In particular, the edge $E$ is orthogonal to its dual facet $F$.
  Finally, the edge $E$ is unbounded if and
  only if its dual facet $F$ is contained in a
  facet of $Q_H$.
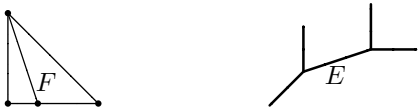
  \begin{center}
    \begin{texdraw}
      \drawdim cm \relunitscale 0.4 \linewd 0.05 \lpatt (1 0) \setgray 0.8 
      \move (0 0) \lvec (3 0) \lvec (0 3) \lvec (0 0) 
      \move (0 0) \fcir f:0 r:0.1 
      \move (0 3) \fcir f:0 r:0.1 
      \move (1 0) \fcir f:0 r:0.1 
      \move (3 0) \fcir f:0 r:0.1 
      \setgray 0
      \move (0 0) \lvec (0 3) \lvec (3 0) \lvec (0 0)
      \move (0 3) \lvec (1 0)
      \htext (0.9 0.4){$F$}
    \end{texdraw}
    \hspace{2cm}
    \begin{texdraw}
      \drawdim cm \relunitscale 0.3 \linewd 0.1 \lpatt (1 0) \setgray 0
      \move (0 0) \lvec (3 1) \rlvec (0 2) \move (3 1) \rlvec (2 0)
      \move (0 0) \rlvec (0 2) \move (0 0) \rlvec (-1.5 -1.5)
      \htext (1 -0.5) {$E$}
    \end{texdraw}
  \end{center}

  The second case we are interested in are the tropical
  hypersurfaces defined by $\trop(A)$ and $\trop(\Delta)$, i.e. the
  case when the tropical hypersurface is actually the
  codimension-one skeleton of a fan. In this case the corresponding
  fan is just the normal fan of the Newton polytope of the the
  defining polynomial (see \cite[Chap.\ 5]{GKZ}).
\end{remark}

\begin{example}\label{ex-markedsubdivision}
  The marked subdivision below is for example induced by the tropical
  polynomial $\min\{3x,3y,0,x,-1+x+y\}$. 
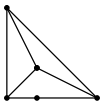
  \begin{center}
    \input{Graphics/exmarksub.pstex_t}
  \end{center}
\end{example}

Let us now restrict our attention to plane tropical cubics, that is, curves defined
by tropical polynomials whose Newton polytope is \emph{contained} in the triangle $Q_3$. Note that
this triangle has only one interior point, so also all possible marked
subdivisions we consider have at most one interior point. 

\begin{definition}\label{def-cycle}
  We say that a plane tropical cubic $\mathcal{C}$ \emph{has a cycle} if the interior
  point $(1,1)$ is visible as the vertex of a marked polytope in its
  dual marked subdivision. If this is the case, the \emph{cycle} of $\mathcal{C}$
  is the union of those bounded edges of $\mathcal{C}$ which are dual to the
  facets of marked polytopes in the marked subdivision which emanate
  from $(1,1)$, and we say that these edges \emph{form the cycle}.
\end{definition}

\begin{example}
  In the picture below, the left plane tropical cubic has a cycle while the
  right one does not, since $(1,1)$ is visible but it is not a vertex of one of
  the marked polytopes in the subdivision.
  \begin{center}
    \input{Graphics/cyclex.pstex_t}
  \end{center}
\end{example}

\begin{definition}\label{def-cyclelength}
  For a bounded edge $E$ of a plane tropical curve with direction vector
  $v(E)$, defined as in Remark \ref{rem-direction} (i.e. $v(E)$ is
  orthogonal to its dual facet in the marked subdivision and of
  the same Euclidean length as this facet) we
  define the \emph{lattice length} $l(E)=\frac{||E||}{||v(E)||}$ to be the Euclidean length of
  $E$ divided by the Euclidean length of $v(E)$.  

  For a plane tropical cubic with cycle, we define its \emph{cycle length} to
  be the sum of the lattice lengths of the edges which form the
  cycle. If the plane tropical cubic has no cycle we say its length is
  zero.
  This defines a function ``cycle length''
  \begin{displaymath}
    \cyclelength:\R^{\ca_3}\rightarrow\R:u=(u_\omega\;|\;\omega\in\ca_3)
    \mapsto\cyclelength(u)=\mbox{``cycle length of $\mathcal{C}_H$''}
  \end{displaymath}
  associating to every plane tropical cubic polynomial
  $H=\min_{\omega\in\ca_3}\{u_\omega+v\cdot\omega\}$ the cycle length of 
  the corresponding plane tropical cubic $\mathcal{C}_H$. 
\end{definition}
\begin{example}
  The following picture shows a plane tropical cubic with cycle length $\frac{9}{2}$.
  \begin{center}
    \input{Graphics/cyclelength.pstex_t}
  \end{center}
\end{example}

\begin{definition}
  Assume a plane tropical cubic has no cycle, i.e.\ $(1,1)$ is not the vertex
  of a marked polytope in the corresponding marked subdivision, but
  it is contained in a facet of (necessarily) two such marked
  polytopes, say $(Q_1,\ca_1)$ and $(Q_2,\ca_2)$.
  They are dual to two vertices $V_1$
  and $V_2$ of the tropical curve. We define the \emph{generalized cycle
    length} of such a cubic to be four times the lattice length of the
  edge connecting $V_1$ and $V_2$. 
\end{definition}

This definition is necessary because these tropical curves arise as
limits of curves with a cycle, and the cycle length tends to the
generalized cycle length for such a limit: 
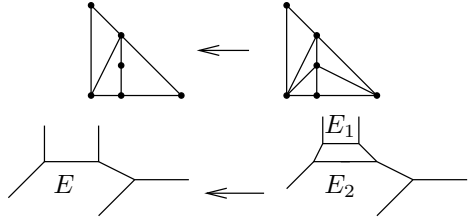
\begin{figure}[h]
  \input{Graphics/cyclelimit.pstex_t}
  \caption{Generalized cycle length}\label{fig-gcl}
\end{figure}
One factor of $2$ is necessary because two edges tend to the same edge
in the picture. The other factor of $2$ appears because the direction
vector of $E$ is of twice the euclidean length than the direction
vector of $E_1$ and $E_2$.

Let us generalize Definition \ref{def-cycle} to plane tropical curves
other than cubics.
\begin{definition}
  Let $\mathcal{C}$ be a plane tropical curve with Newton polytope $Q$ and with dual marked
  subdivision $\{(Q_i,\ca_i)\;|\;i=1,\ldots,l\}$. Suppose
  that $\tilde{\omega}\in\Int(Q)\cap\Z^2$ is in the interior of $Q$ and that the
  $(Q_i,\ca_i)$ are ordered such that $\tilde{\omega}$ is a vertex of $Q_i$
  for $i=1,\ldots,k$ and it is not contained in $Q_i$ for
  $i=k+1,\ldots,l$ (see Figure \ref{fig-cycle}). 
  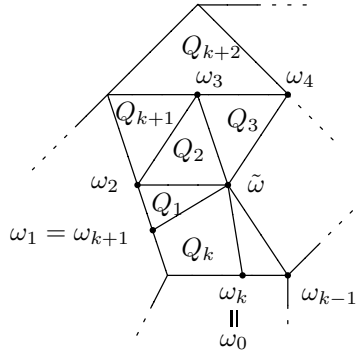
\begin{figure}[h]
    \centering
    \begin{texdraw}
      \drawdim cm  \relunitscale 0.4 \linewd 0.05  \lpatt (1 0) \setgray 0
      \move (0 0) \rlvec (2 3) \rlvec (3 0)
      \lpatt (0.1 0.4) \rlvec (1.5 -1.5) \lpatt(1 0)
      \move (0 0) \rlvec (1 -3) \rlvec (4 0) \rlvec (1 1)
      \lpatt (0.1 0.4) \rlvec (1.3 1.3) \lpatt(1 0)
      \move (3 0) \lvec (0 0) \fcir f:0 r:0.1 \htext (-1.5 -0.2) {$\omega_2$}
      \move (3 0) \lvec (0.5 -1.5) \fcir f:0 r:0.1  \htext (-4.2 -2) {$\omega_1=\omega_{k+1}$}
      \move (3 0) \lvec (3.5 -3) \fcir f:0 r:0.1 
      \htext (2.5 -5.7) {$\xymatrix@R0.2cm{\omega_k\ar@{=}[d]\\\omega_0}$}
      \move (3 0) \lvec (5 -3) \fcir f:0 r:0.1  \htext (5.5 -4) {$\omega_{k-1}$}
      \move (3 0) \lvec (2 3) \fcir f:0 r:0.1 \htext (2 3.3) {$\omega_3$}
      \move (3 0) \lvec (5 3) \fcir f:0 r:0.1 \htext (5 3.3) {$\omega_4$}
      \move (3 0) \fcir f:0 r:0.1 \htext (3.7 -0.2) {$\tilde{\omega}$}
      \htext (1.5 -2.5) {$Q_k$}
      \htext (0.5 -1) {$Q_1$}
      \htext (1.2 0.8) {$Q_2$}
      \htext (3 1.8) {$Q_3$}
      \move (0 0) \rlvec (-1 3) \lvec (2 3) \htext (-0.7 2) {$Q_{k+1}$}
      \move (-1 3) \rlvec (3 3) \rlvec (3 -3) \htext (1.5 4.2) {$Q_{k+2}$}
      \move (2 6) \rlvec (2 0) \lpatt (0.1 0.4) \rlvec (2 0) 
      \move (-1 3) \lpatt (1 0) \rlvec (-1 -1) \lpatt (0.1 0.4)
      \rlvec (-1.5 -1.5)
      \move (1 -3) \lpatt (1 0) \rlvec (-0.5 -1) \lpatt (0.1 0.4)
      \rlvec (-0.5 -1)
      \move (5 -3) \lpatt (1 0) \rlvec (0 -1) \lpatt (0.1 0.4)
      \rlvec (0 -1)
      
    \end{texdraw}    
    \caption{Marked subdivision determining a cycle}
    \label{fig-cycle}
  \end{figure}
  We then say that $\tilde{\omega}$ \emph{determines a
    cycle} of $\mathcal{C}$, namely the union of the edges of $\mathcal{C}$ dual to the
  facets emanating from $\tilde{\omega}$, and we say that these edges 
  \emph{form the cycle} determined by $\tilde{\omega}$. The length of
  this cycle is defined as in  Definition \ref{def-cyclelength}.
\end{definition}

\begin{lemma}\label{lem-cyclelengthlinear}
  Let $(Q,\mathcal{A})$  be a marked polytope in $\R^2$
  with a regular marked subdivision $\{(Q_i,\ca_i)\;|\;i=1,\ldots,l\}$
  and suppose that $\tilde{\omega}\in\Int(Q)\cap\Z^2$ in the interior of $Q$
  is a vertex of $Q_i$ for $i=1,\ldots,k$ and it is not contained in
  $Q_i$ for $i=k+1,\ldots,l$. 

  If $\psi:\mathcal{A}\rightarrow \R$ is a function defining the
  subdivision, then $\tilde{\omega}$ determines a cycle in the plane 
  tropical curve $\mathcal{C}$ given by the tropical polynomial
  \begin{displaymath}
    \min\{\psi(\omega)+v\cdot\omega\;|\;\omega\in\mathcal{A}\}
  \end{displaymath}
  and, using the notation in Figure \ref{fig-cycle}, its length is
  \begin{displaymath}
    \sum_{j=1}^k\big(\psi(\tilde{\omega})-\psi(\omega_j)\big)\cdot
    \frac{D_{j-1,j}+D_{j,j+1}+D_{j+1,j-1}}{D_{j-1,j}\cdot D_{j,j+1}}
  \end{displaymath}
  where $D_{i,j}=\det(w_i, w_j)$ with $w_i=\omega_i-\tilde{\omega}$ and $w_j=\omega_j-\tilde{\omega}$.
\end{lemma}
\begin{proof}
  By definition $\tilde{\omega}$ determines a cycle. It remains to prove the
  statement on its length.
  
  For this we consider the convex polytope $Q_j$ having $\omega_{j+1}$,
  $\tilde{\omega}$ and $\omega_{j}$ as neighboring vertices:
  \smallskip
  
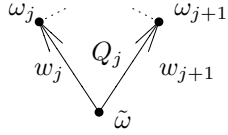
\begin{center}
    \begin{texdraw}
      \drawdim cm  \relunitscale 0.2 \linewd 0.05  \lpatt (1 0)
      \setgray 0  \arrowheadtype t:V
      \move (0 0) \ravec (-4 6)  \fcir f:0 r:0.3
      \move (0 0) \ravec (4 6)   \fcir f:0 r:0.3  
      \move (0 0) \fcir f:0 r:0.3
      \htext (1 -1) {$\tilde{\omega}$}
      \htext (-6 6) {$\omega_j$}
      \htext (5 6) {$\omega_{j+1}$}
      \htext (-4.3 2) {$w_j$}
      \htext (4 2) {$w_{j+1}$}
      \htext (-0.5 3) {$Q_j$}
      \lpatt (0.1 0.4) \move (-4 6) \rlvec (2 1)
      \move (4 6) \rlvec (-2 1)
    \end{texdraw}
  \end{center}
  The vertex $v_j=(v_{j,1},v_{j,2})$ of $\mathcal{C}$ dual to $Q_j$ is given by
  the system of linear equations
  \begin{displaymath}
    \omega_j\cdot v_j  +u_j=
    \omega_{j+1}\cdot v_j +u_{j+1}=
    \tilde{\omega}\cdot v_j +u,
  \end{displaymath}
  where $u_j=\psi(\omega_j)$, $u_{j+1}=\psi(\omega_{j+1})$ and
  $u=\psi(\tilde{\omega})$.
  This system can be rewritten as
  \begin{displaymath}
    \begin{pmatrix}
      w_j^t\\ w_{j+1}^t
    \end{pmatrix}
    \cdot v_j=
    \begin{pmatrix}
      u-u_j\\u-u_{j+1}
    \end{pmatrix}.
  \end{displaymath}
  Since $\omega_{j+1}$, $\tilde{\omega}$ and $\omega_j$ are neighboring
  vertices of the polytope $Q_j$ the vectors $w_j$ and $w_{j+1}$ are
  linearly independent and we may apply Cramer's Rule to find
  \begin{equation}\label{eq:latticlength:0}
    v_{j,1}=\frac{\DET{cc}{u-u_j & w_{j,2}\\ u-u_{j+1} & w_{j+1,2}}}{D_{j,j+1}}
    \;\;\; \mbox{ and }\;\;\;
    v_{j,2}=\frac{\DET{cc}{w_{j,1} & u-u_j  \\ w_{j+1,1} & u-u_{j+1}}}{D_{j,j+1}}.
  \end{equation}
  The lattice length of the edge from $v_{j-1}$ to $v_j$ is the real
  number $\lambda_j\in\R$ such that $(v_j-v_{j-1})=\lambda_j\cdot
  w_j^\perp$, where where $w_j^\perp=(-w_{j,2},w_{j,1})$ 
  is perpendicular to $w_j$. Thus
  \begin{equation}\label{eq:latticelength:1}
    \lambda_j=
    \frac{ (v_j-v_{j-1})\cdot w_j^\perp}{ w_j^\perp\cdot w_j^\perp}=
    \frac{(v_j-v_{j-1})\cdot w_j^\perp}{ w_j\cdot w_j}.
  \end{equation}
  In order to understand the right hand side of this equation better we need
  the following observation. The last row of the matrix
  \begin{displaymath}
    M=\begin{pmatrix} w_{j-1,1} & w_{j,1} & w_{j+1,1}\\
    w_{j-1,2} & w_{j,2} & w_{j+1,2}\\
      w_{j-1}\cdot w_j &
      w_j\cdot w_j &
      w_{j+1}\cdot w_j
    \end{pmatrix}
  \end{displaymath}
  is a linear combination of the first two, and thus the determinant
  of $M$ is zero. Developing the determinant by the last row we get
  \begin{displaymath}
    0=\det(M)=
    w_{j-1}\cdot w_j \cdot D_{j,j+1} -
    w_j\cdot w_j \cdot D_{j-1,j+1} +
    w_{j+1}\cdot w_j \cdot D_{j-1,j},
  \end{displaymath}
  or equivalently
  \begin{displaymath}
    \frac{D_{j+1,j-1}}{D_{j-1,j}\cdot D_{j,j+1}}=
    -\frac{D_{j-1,j+1}}{D_{j-1,j}\cdot D_{j,j+1}}=
    -\frac{ w_{j-1}\cdot w_j}{ w_j\cdot w_j\cdot D_{j-1,j}} 
    -\frac{ w_{j+1}\cdot w_j}{ w_j\cdot w_j\cdot D_{j,j+1}}.
  \end{displaymath}
  Expanding the right hand side of \eqref{eq:latticelength:1} using
  \eqref{eq:latticlength:0} and
  plugging in this last equality we get
  \begin{align*}
    \lambda_j=&
    \frac{u-u_{j-1}}{D_{j-1,j}}+
    \frac{u-u_{j+1}}{D_{j,j+1}}-
    (u-u_j)\cdot\left(
      \frac{ w_j\cdot w_{j+1}}{ w_j\cdot w_j\cdot D_{j,j+1}}+
      \frac{ w_{j-1}\cdot w_j}{ w_j\cdot w_j\cdot D_{j-1,j}}\right)      \\
    =&
    \frac{u-u_{j-1}}{D_{j-1,j}}+\frac{u-u_{j+1}}{D_{j,j+1}}
    +\frac{(u-u_j)\cdot D_{j+1,j-1}}{D_{j-1,j}\cdot D_{j,j+1}}.
  \end{align*}
  The lattice length of the cycle of $\mathcal{C}$ is then given by adding the
  $\lambda_j$, i.e. it is
  \begin{align*}
    \lambda_1+\ldots+\lambda_k=&
    \sum_{j=1}^k
    \frac{u-u_{j-1}}{D_{j-1,j}}+\frac{u-u_{j+1}}{D_{j,j+1}}+
    \frac{(u-u_j)\cdot D_{j+1,j-1}}{D_{j-1,j}\cdot D_{j,j+1}}\\
    =&
    \sum_{j=1}^k (u-u_j)\cdot
    \left(\frac{D_{j-1,j}+D_{j,j+1}+D_{j+1,j-1}}{D_{j-1,j}\cdot D_{j,j+1}}\right).
  \end{align*}
\end{proof}

\begin{remark}\label{rem-cyclelengthlinear}
  An immediate consequence of Lemma \ref{lem-cyclelengthlinear} is that
  the function ``cycle length'', $\cyclelength$, from Definition \ref{def-cyclelength} is
  linear on each cone of the secondary fan of $\ca_3$. 
\end{remark}


\section{The main theorem}\label{sec-thm}

\begin{theorem}\label{thm-main}
  Let $\mathcal{C}$ be a plane tropical cubic  given by the tropical
  polynomial 
  \begin{displaymath}
    \min_{(i,j)\in\ca_3}\{u_{ij}+ix+jy\}
  \end{displaymath}
  and assume that $\mathcal{C}$ has a cycle.

  Then the negative of the generic valuation of the $j$-invariant at
  $u=(u_{ij})_{(i,j)\in\ca_3}$ 
  is equal to the cycle length of $\mathcal{C}$, i.e.\
  \begin{displaymath}
    -\val_u(j)=\cyclelength(u).
  \end{displaymath}

  Furthermore, if the marked subdivision dual to $\mathcal{C}$
  corresponds to a
  top-dimensional cone of the secondary fan of $\ca_3$ (that is, if it is a
  triangulation), then $\val_u(j)=\val(j(f))$ where
  $f=\sum_{(i,j)\in\ca_3}a_{ij}x^iy^j$ is any elliptic curve over $\K$
  with coefficients $a_{ij}$ satisfying $\val(a_{ij})=u_{ij}$. 
\end{theorem}

There are two main parts of the proof: the first part is to compare
certain ``domains of linearity'' in $\R^{\ca_3}$ of the two piece-wise
linear functions ``cycle length'', $\cyclelength$, and ``generic valuation of
$j$'', $\val_{\cdot}(j)$, and
the second part is to compare the two linear functions on each
domain. 

The proof uses many results that will be proved in the following sections.

\begin{proof}[Proof of Theorem \ref{thm-main}:]
  Note that our claim only involves curves $\mathcal{C}$ which have a cycle
  or, equivalently, where in the dual subdivision the point $(1,1)$ is a
  vertex of a marked polytope. Therefore we may replace $\R^{\ca_3}$ as
  domain of definition of $\cyclelength$ and $\val_\cdot(j)$ by the
  union $U$ of
  those cones of the secondary fan of $\ca_3$ where the corresponding
  marked subdivision contains $(1,1)$ as a vertex of a marked polytope.
  The coordinates on $U$ are given by
  $u=\big(u_{ij}\;|\;(i,j)\in\ca_3\big)$ and the canonical basis vector
  $e_{kl}=\big(\delta_{ik}\cdot \delta_{jl}\;|\;(i,j)\in\ca_3\big)$ has a
  one in position $kl$ and zeros elsewhere.

  From Lemma \ref{lem-num} we know that $U$
  is contained in a single cone of the Gr\"obner fan of $A$, namely the
  one dual to the vertex $12\cdot e_{11}$ of the Newton polytope of $A$. 
  Hence the generic valuation of $A$ is linear on $U$, namely 
  \begin{displaymath}
    U\rightarrow\R:u\mapsto\val_u(A)=12\cdot u_{11}.  
  \end{displaymath}
  Thus, if we want to divide $U$ into cones on which $\val_{\cdot}(j)$
  is linear, it suffices to consider $u\mapsto\val_u(\Delta)$, and we
  know already that the latter is linear on cones of the Gr\"obner fan
  of $\Delta$ by Lemma \ref{lem-genvallinear}. Thus so is
  $\val_{\cdot}(j)$ restricted to $U$, 
  and by Lemma 
  \ref{lem-cyclinondelta} and Remark \ref{rem-deltaequivalence} 
  the function $\cyclelength$ is so as well. Moreover, Remark
  \ref{rem-deltaequivalence} tells us that $U$ is indeed a union of
  cones of the Gr\"obner fan of $\Delta$, and each such cone is a
  union of certain $\Delta$-equivalent cones of the secondary fan of $\ca_3$.

  Hence to prove that the two functions $\val_{\cdot}(j)$ and $\cyclelength$ 
  coincide it is
  enough to compare the linear functions on each cone of the Gr\"obner
  fan of $\Delta$ contained in $U$. To do this, we
  use Theorem 11.3.2 of \cite{GKZ} which enables us to compute the
  assignment rule for the linear function $u\mapsto\val_u(\Delta)$ on each such cone, say $D$,
  given a (top-dimensional) marked subdivision $T$ whose corresponding cone in
  the secondary fan of $\ca_3$ is contained in $D$. In fact, it provides us
  with a formula to compute the coefficient of $u_{ij}$ for each
  $(i,j)\in \ca_3$. 
  Since we already know that the two functions $u\mapsto\val_u(\Delta)$
  and $\cyclelength$ are linear on $D$, we can for
  our comparison assume that $T$ is the representative of its class with
  as few edges as possible.  
  The coefficient of $u_{ij}$ in the linear function $\cyclelength$ for
  the marked subdivision is given by Lemma \ref{lem-cyclelengthlinear}. 
  To compare the two coefficients, there are some cases to distinguish,
  which is done by Lemma \ref{lem-comparelin}. This proves the first part
  of the theorem.

  Finally, for any point $u$ in the interior of a cone of the Gr\"obner fan of
  $\Delta$, $\val_u(j)=\val(j(f))$ for any polynomial $f=\sum_{(i,j)\in\ca_3}
  a_{ij}x^i y^j$ with $\val(a_{ij})=u_{ij}$ by Lemma
  \ref{lem-genvallinear}. As a point $u$ in the interior of a
  top-dimensional cone of the secondary fan of $\ca_3$ is in the interior of a cone
  of the Gr\"obner fan of $\Delta$, the last statement follows as well. 
\end{proof}

We would like to give an alternative proof of the statement whose
methods we believe to be interesting on their own. 
Here, we will consider smaller domains of linearity, namely the cones
of the secondary fan of $\ca_3$ contained in $U$ (using the notation
from the above proof). 

\begin{proof}[Alternative proof of Theorem \ref{thm-main}:]
  From \cite[Chap.\ 10, Thm.\ 1.2]{GKZ} or alternatively from Lemma
  \ref{discriminant} we conclude that the codimension $1$ cones of the
  Gr\"obner fan of $\Delta$ do no meet the interior of any
  top-dimensional cone of the secondary fan of $\ca_3$. Thus an open
  top-dimensional cone of the secondary fan of $\ca_3$ is completely
  contained in some top-dimensional cone of the Gr\"obner fan of
  $\Delta$, and using Lemma \ref{lem-genvallinear} we conclude that
  $u\mapsto\val_u(\Delta)$ is linear on each top-dimensional cone of the
  secondary fan of $\ca_3$. Using Lemma \ref{lem-num} we can see that
  $u\mapsto\val_u(A)$ is linear on each cone of the secondary fan of
  $\ca_3$ corresponding to a subdivision for which the interior point is visible.
    By \ref{lem-cyclelengthlinear} we know that the function
  $\cyclelength$ is linear on a cone of the secondary fan of $\ca_3$, too. To
  show that the two functions agree, we only have to show that they
  agree on the rays of each cone of the secondary fan of $\ca_3$ in question. In
  Proposition \ref{prop-rays} we classify the rays. Then a computation
  for each ray shows that the two functions agree. We computed this
  using the procedure \texttt{raysC} in the library
  \texttt{jinvariant.lib} available via the URL
  \begin{center}
    http://www.mathematik.uni-kl.de/\textasciitilde keilen/en/jinvariant.html.
  \end{center}
  Note that we have to compare with the generalized cycle length,
  because a point on a ray is a limit of points corresponding to curves
  with a cycle. Since the rays of the secondary fan of $\ca_3$ in question are not
  necessarily contained in the interior of a top-dimensional cone of the Gr\"obner
  fan of $A$ and $\Delta$, we have to use the generic
  valuation. However, for a point $u$ in the interior of a
  top-dimensional cone of the secondary fan of $\ca_3$ we know that
  $\val_u(j)=\val(j(f))$ for any polynomial
  $f=\sum_{(i,j)\in\ca_3} a_{ij}x^i y^j$ with $\val(a_{ij})=u_{ij}$ by
  Lemma \ref{lem-genvallinear}.  
\end{proof}

\begin{corollary}\label{cor-nocycle}
  Let $f=\sum_{(i,j)\in\ca_3} a_{ij}x^i y^j$ define a smooth elliptic
  curve over $\K$ such that the valuation of its $j$-invariant is
  positive. Then its tropicalization does not have a cycle. 
\end{corollary}

\begin{remark}
  Let $C$ be a smooth elliptic curve over $\K$. It is obvious that the
  tropicalization of $C$ depends on the embedding into the projective
  plane that we choose. One might ask if the cycle length does not
  depend on the embedding though, as it should take the role of a
  tropical $j$-invariant. This is not true however. For example, each
  curve can be put into Weierstrass-form without $xy$-term and its tropicalization does
  not have a cycle, because the interior point is not part of the marked
  subdivision. 
  Also, an embedding might be such that the valuations of the
  coefficients lie in a cone of the Gr\"obner fan of $\Delta$ of higher
  codimension, and the coefficients might be such that the lowest order
  terms of $\Delta(f)$ cancel. But as the cycle length is equal to the
  generic valuation of $j$ and the generic valuation of $j$ is smaller than the
  valuation of the $j$-invariant in this case, the cycle length would
  not reflect the $j$-invariant. 

  As an example, we consider the following family of curves over $\K$
  with the same $j$-invariant (hence you could also say: a family of
  different embeddings for one curve). We choose a given curve and apply
  the coordinate change (in affine coordinates) $(x,y)\mapsto (x+k,y)$, where
  $k\in\K$. Since this is an 
  isomorphism, the curves over $\K$ have the same $j$-invariant for any
  choice of $k$. In particular, the valuation of the $j$-invariant is
  the same for any $k$. 
  One might hope that at least for a general choice of $k$ the cycle
  length of the tropicalization is equal to the valuation of the
  $j$-invariant. But even this is not true, as the following example
  demonstrates.  
  Let us take a subset of the family given by the coordinate changes
  $(x,y)\mapsto (x+t^b,y)$, where $b\in \Q$. For our example we will see that
  infinitely many choices of $b$ lead to the expected cycle length, but
  also infinitely many lead to the ``wrong'' cycle length. This shows
  that it is not true that a general choice of $k$ (in the sense of
  Zariski topology) leads to the expected
  cycle length. 
  As example we choose 
  \begin{align*}f=& c_{00}\cdot t+c_{10}\cdot
    t^{100}\cdot x+c_{20}\cdot t^{100}\cdot x^2+c_{30}\cdot t\cdot
    x^3+c_{01}\cdot t\cdot y\\&+c_{11}\cdot x\cdot y+c_{21}\cdot
    t^{100}\cdot x^2\cdot y+c_{02}\cdot t^3\cdot y^2+c_{12}\cdot t\cdot
    x\cdot y^2+c_{03}\cdot t^7\cdot y^3,
  \end{align*} 
  where the $c_{ij}\in\C$ are general.  
  By general we mean that after applying the coordinate change every
  coefficient has the expected valuation and nothing cancels. 
  The valuation of the $j$-invariant is $5$. The cycle length of the
  tropicalization is $5$, too, as expected. 

  
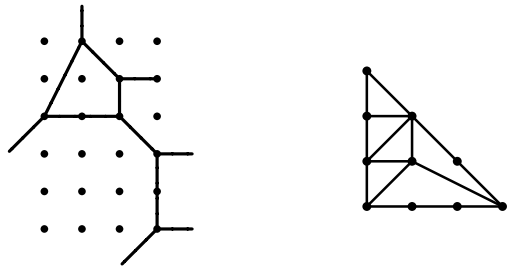
\begin{center}
    \begin{texdraw}
       \drawdim cm  \relunitscale 0.2 \arrowheadtype t:V
       \linewd 0.2  \lpatt (1 0)      
       \relunitscale 2.5
       \move (-1 0) \lvec (0 2)
       \move (-1 0) \lvec (1 0)
       \move (-1 0) \rlvec (-0.93 -0.93)
       \move (0 2) \lvec (1 1)
       \move (0 2) \rlvec (0 0.93)
       \move (2 -3) \lvec (2 -1)
       \move (2 -3) \rlvec (0.93 0)
       \move (2 -3) \rlvec (-0.93 -0.93)
       \move (1 0) \lvec (2 -1)
       \move (1 0) \lvec (1 1)
       \move (2 -1) \rlvec (0.93 0)
       \move (1 1) \rlvec (0.93 0)

        \move (-1 -3) \fcir f:0 r:0.1
        \move (-1 -2) \fcir f:0 r:0.1
        \move (-1 -1) \fcir f:0 r:0.1
        \move (-1 0) \fcir f:0 r:0.1
        \move (-1 1) \fcir f:0 r:0.1
        \move (-1 2) \fcir f:0 r:0.1
        \move (0 -3) \fcir f:0 r:0.1
        \move (0 -2) \fcir f:0 r:0.1
        \move (0 -1) \fcir f:0 r:0.1
        \move (0 0) \fcir f:0 r:0.1
        \move (0 1) \fcir f:0 r:0.1
        \move (0 2) \fcir f:0 r:0.1
        \move (1 -3) \fcir f:0 r:0.1
        \move (1 -2) \fcir f:0 r:0.1
        \move (1 -1) \fcir f:0 r:0.1
        \move (1 0) \fcir f:0 r:0.1
        \move (1 1) \fcir f:0 r:0.1
        \move (1 2) \fcir f:0 r:0.1
        \move (2 -3) \fcir f:0 r:0.1
        \move (2 -2) \fcir f:0 r:0.1
        \move (2 -1) \fcir f:0 r:0.1
        \move (2 0) \fcir f:0 r:0.1
        \move (2 1) \fcir f:0 r:0.1
        \move (2 2) \fcir f:0 r:0.1
    \end{texdraw}
    \hspace*{2cm}
    \begin{texdraw}
       \drawdim cm  \relunitscale 0.6
       \linewd 0.05
        \move (3 0)        
        \lvec (0 0)
        \move (0 0)        
        \lvec (0 3)
        \move (0 3)        
        \lvec (3 0)    
        \move (3 0)        
        \lvec (1 1)
        \move (1 1)        
        \lvec (1 2)
        \move (0 0)        
        \lvec (1 1)
        \move (1 2)        
        \lvec (0 2)
        \move (0 1)        
        \lvec (1 2)
        \move (1 1)        
        \lvec (0 1)
        \move (0 0) \fcir f:0 r:0.1
        \move (0 1) \fcir f:0 r:0.1
        \move (0 2) \fcir f:0 r:0.1
        \move (0 3) \fcir f:0 r:0.1
        \move (1 0) \fcir f:0 r:0.1
        \move (1 1) \fcir f:0 r:0.1
        \move (1 2) \fcir f:0 r:0.1
        \move (2 0) \fcir f:0 r:0.1
        \move (2 1) \fcir f:0 r:0.1
        \move (3 0) \fcir f:0 r:0.1
        \move (0 -1.2) \fcir f:1 r:0.1
    \end{texdraw}
  \end{center}

  Let us check what happens to the valuations $u_{ij}$ of the
  coefficients $c_{ij}t^{u_ij}$ when we apply the coordinate change. In
  general (i.e.\ if no cancellation happens) the valuations are as
  follows: 

  \hspace*{0.5cm}
  \input{Graphics/extrans.pstex_t}

  If we choose $b$ very small, then all points but the point $(0,3)$ lie
  on the same plane and thus the subdivision is as follows: 

  \begin{center}
    \input{Graphics/raysec.pstex_t}
  \end{center}
  In fact, it corresponds to a ray of the secondary fan of $\ca_3$.

  For our special choice of $f$ and the coefficients $u_{ij}$ as above,
  this marked subdivision is only reached when $b\leq -1$. Starting from the
  other end we see that the tropical curve stays unchanged as long as $b\geq 2$.  
  In particular, the cycle length is as expected for all those $b$.
  For $b$ in the interval $1< b\leq 2$ the dual marked subdivision is
  not changed, but the position of two vertices of the tropical curve
  changes. For example for $b=\frac{3}{2}$, the tropical curve is: 

  
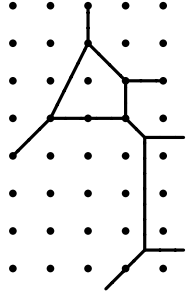
\begin{center}
    \begin{texdraw}
       \drawdim cm  \relunitscale 0.2 \arrowheadtype t:V
       \linewd 0.2  \lpatt (1 0)      
       \relunitscale 2.5
       \move (-1 0) \lvec (0 2)
       \move (-1 0) \lvec (1 0)
       \move (-1 0) \rlvec (-1.03 -1.03)
       \move (0 2) \lvec (1 1)
       \move (0 2) \rlvec (0 1.03)
       \move (1.5 -3.5) \lvec (1.5 -0.5)
       \move (1.5 -3.5) \rlvec (1.03 0)
       \move (1.5 -3.5) \rlvec (-1.03 -1.03)
       \move (1 0) \lvec (1.5 -0.5)
       \move (1 0) \lvec (1 1)
       \move (1.5 -0.5) \rlvec (1.03 0)
       \move (1 1) \rlvec (1.03 0)

        \move (-2 -4) \fcir f:0 r:0.1
        \move (-2 -3) \fcir f:0 r:0.1
        \move (-2 -2) \fcir f:0 r:0.1
        \move (-2 -1) \fcir f:0 r:0.1
        \move (-2 0) \fcir f:0 r:0.1
        \move (-2 1) \fcir f:0 r:0.1
        \move (-2 2) \fcir f:0 r:0.1
        \move (-2 3) \fcir f:0 r:0.1
        \move (-1 -4) \fcir f:0 r:0.1
        \move (-1 -3) \fcir f:0 r:0.1
        \move (-1 -2) \fcir f:0 r:0.1
        \move (-1 -1) \fcir f:0 r:0.1
        \move (-1 0) \fcir f:0 r:0.1
        \move (-1 1) \fcir f:0 r:0.1
        \move (-1 2) \fcir f:0 r:0.1
        \move (-1 3) \fcir f:0 r:0.1
        \move (0 -4) \fcir f:0 r:0.1
        \move (0 -3) \fcir f:0 r:0.1
        \move (0 -2) \fcir f:0 r:0.1
        \move (0 -1) \fcir f:0 r:0.1
        \move (0 0) \fcir f:0 r:0.1
        \move (0 1) \fcir f:0 r:0.1
        \move (0 2) \fcir f:0 r:0.1
        \move (0 3) \fcir f:0 r:0.1
        \move (1 -4) \fcir f:0 r:0.1
        \move (1 -3) \fcir f:0 r:0.1
        \move (1 -2) \fcir f:0 r:0.1
        \move (1 -1) \fcir f:0 r:0.1
        \move (1 0) \fcir f:0 r:0.1
        \move (1 1) \fcir f:0 r:0.1
        \move (1 2) \fcir f:0 r:0.1
        \move (1 3) \fcir f:0 r:0.1
        \move (2 -4) \fcir f:0 r:0.1
        \move (2 -3) \fcir f:0 r:0.1
        \move (2 -2) \fcir f:0 r:0.1
        \move (2 -1) \fcir f:0 r:0.1
        \move (2 0) \fcir f:0 r:0.1
        \move (2 1) \fcir f:0 r:0.1
        \move (2 2) \fcir f:0 r:0.1
        \move (2 3) \fcir f:0 r:0.1
    \end{texdraw}
  \end{center}

  For $b=1$ the marked subdivision changes to

  
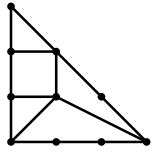
\begin{center}    
    \begin{texdraw}
       \drawdim cm  \relunitscale 0.6
       \linewd 0.05
        \move (3 0)        
        \lvec (0 0)
        \move (0 0)        
        \lvec (0 3)
        \move (0 3)        
        \lvec (3 0)

        \move (3 0)        
        \lvec (1 1)
        \move (1 1)        
        \lvec (1 2)
        \move (0 0)        
        \lvec (1 1)
        \move (1 2)        
        \lvec (0 2)
        \move (1 1)        
        \lvec (0 1)
        \move (0 0) \fcir f:0 r:0.08
        \move (0 1) \fcir f:0 r:0.08
        \move (0 2) \fcir f:0 r:0.08
        \move (0 3) \fcir f:0 r:0.08
        \move (1 0) \fcir f:0 r:0.08
        \move (1 1) \fcir f:0 r:0.08
        \move (1 2) \fcir f:0 r:0.08
        \move (2 0) \fcir f:0 r:0.08
        \move (2 1) \fcir f:0 r:0.08
        \move (3 0) \fcir f:0 r:0.08
    \end{texdraw}     
  \end{center}
  and it remains like this for $0\leq b\leq 1$, while
  the cycle length decreases; e.g. for $b=\frac{2}{3}$ respectively
  $b=\frac{1}{3}$ we get:
  
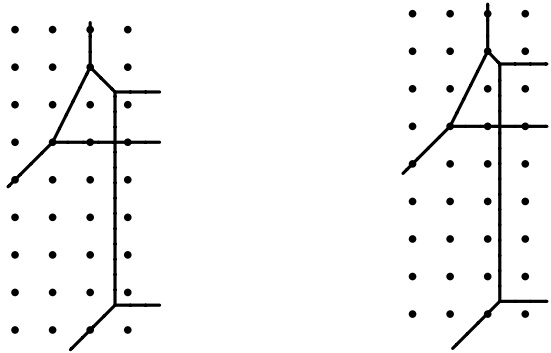
\begin{center}
    \begin{texdraw}
       \drawdim cm  \relunitscale 0.2 \arrowheadtype t:V
       \linewd 0.2  \lpatt (1 0)       
       \relunitscale 2.5
       \move (-1 1) \lvec (0 3)
       \move (-1 1) \lvec (0.66 1)
       \move (-1 1) \rlvec (-1.18 -1.18)
       \move (0 3) \lvec (0.66 2.33)
       \move (0 3) \rlvec (0 1.18)
       \move (0.66 -3.33) \lvec (0.66 1)
       \move (0.66 -3.33) \rlvec (1.18 0)
       \move (0.66 -3.33) \rlvec (-1.18 -1.18)
       \move (0.66 1) \lvec (0.66 2.33)
       \move (0.66 1) \rlvec (1.18 0)
       \move (0.66 2.33) \rlvec (1.18 0)

        \move (-2 -4) \fcir f:0 r:0.1
        \move (-2 -3) \fcir f:0 r:0.1
        \move (-2 -2) \fcir f:0 r:0.1
        \move (-2 -1) \fcir f:0 r:0.1
        \move (-2 0) \fcir f:0 r:0.1
        \move (-2 1) \fcir f:0 r:0.1
        \move (-2 2) \fcir f:0 r:0.1
        \move (-2 3) \fcir f:0 r:0.1
        \move (-2 4) \fcir f:0 r:0.1
        \move (-1 -4) \fcir f:0 r:0.1
        \move (-1 -3) \fcir f:0 r:0.1
        \move (-1 -2) \fcir f:0 r:0.1
        \move (-1 -1) \fcir f:0 r:0.1
        \move (-1 0) \fcir f:0 r:0.1
        \move (-1 1) \fcir f:0 r:0.1
        \move (-1 2) \fcir f:0 r:0.1
        \move (-1 3) \fcir f:0 r:0.1
        \move (-1 4) \fcir f:0 r:0.1
        \move (0 -4) \fcir f:0 r:0.1
        \move (0 -3) \fcir f:0 r:0.1
        \move (0 -2) \fcir f:0 r:0.1
        \move (0 -1) \fcir f:0 r:0.1
        \move (0 0) \fcir f:0 r:0.1
        \move (0 1) \fcir f:0 r:0.1
        \move (0 2) \fcir f:0 r:0.1
        \move (0 3) \fcir f:0 r:0.1
        \move (0 4) \fcir f:0 r:0.1
        \move (1 -4) \fcir f:0 r:0.1
        \move (1 -3) \fcir f:0 r:0.1
        \move (1 -2) \fcir f:0 r:0.1
        \move (1 -1) \fcir f:0 r:0.1
        \move (1 0) \fcir f:0 r:0.1
        \move (1 1) \fcir f:0 r:0.1
        \move (1 2) \fcir f:0 r:0.1
        \move (1 3) \fcir f:0 r:0.1
        \move (1 4) \fcir f:0 r:0.1
    \end{texdraw}
    \hspace{3cm}
    \begin{texdraw}
       \drawdim cm  \relunitscale 0.2 \arrowheadtype t:V
       \linewd 0.2  \lpatt (1 0)       
       \relunitscale 2.5
       \move (-1 1) \lvec (0 3)
       \move (-1 1) \lvec (0.33 1)
       \move (-1 1) \rlvec (-1.25 -1.25)
       \move (0 3) \lvec (0.33 2.66)
       \move (0 3) \rlvec (0 1.25)
       \move (0.33 -3.66) \lvec (0.33 1)
       \move (0.33 -3.66) \rlvec (1.25 0)
       \move (0.33 -3.66) \rlvec (-1.25 -1.25)
       \move (0.33 1) \lvec (0.33 2.66)
       \move (0.33 1) \rlvec (1.25 0)
       \move (0.33 2.66) \rlvec (1.25 0)

        \move (-2 -4) \fcir f:0 r:0.1
        \move (-2 -3) \fcir f:0 r:0.1
        \move (-2 -2) \fcir f:0 r:0.1
        \move (-2 -1) \fcir f:0 r:0.1
        \move (-2 0) \fcir f:0 r:0.1
        \move (-2 1) \fcir f:0 r:0.1
        \move (-2 2) \fcir f:0 r:0.1
        \move (-2 3) \fcir f:0 r:0.1
        \move (-2 4) \fcir f:0 r:0.1
        \move (-1 -4) \fcir f:0 r:0.1
        \move (-1 -3) \fcir f:0 r:0.1
        \move (-1 -2) \fcir f:0 r:0.1
        \move (-1 -1) \fcir f:0 r:0.1
        \move (-1 0) \fcir f:0 r:0.1
        \move (-1 1) \fcir f:0 r:0.1
        \move (-1 2) \fcir f:0 r:0.1
        \move (-1 3) \fcir f:0 r:0.1
        \move (-1 4) \fcir f:0 r:0.1
        \move (0 -4) \fcir f:0 r:0.1
        \move (0 -3) \fcir f:0 r:0.1
        \move (0 -2) \fcir f:0 r:0.1
        \move (0 -1) \fcir f:0 r:0.1
        \move (0 0) \fcir f:0 r:0.1
        \move (0 1) \fcir f:0 r:0.1
        \move (0 2) \fcir f:0 r:0.1
        \move (0 3) \fcir f:0 r:0.1
        \move (0 4) \fcir f:0 r:0.1
        \move (1 -4) \fcir f:0 r:0.1
        \move (1 -3) \fcir f:0 r:0.1
        \move (1 -2) \fcir f:0 r:0.1
        \move (1 -1) \fcir f:0 r:0.1
        \move (1 0) \fcir f:0 r:0.1
        \move (1 1) \fcir f:0 r:0.1
        \move (1 2) \fcir f:0 r:0.1
        \move (1 3) \fcir f:0 r:0.1
        \move (1 4) \fcir f:0 r:0.1
    \end{texdraw}
  \end{center}
  This happens because the cycle length is equal to the generic
  valuation of $j$. The latter is not equal to the valuation of the
  $j$-invariant if and only if the $t$-initial form $\tini_u(\Delta)$
  cancels when plugging in the leading terms of the Puiseux series
  coefficients.  
  The $t$-initial form of $\Delta$ corresponding to this marked
  subdivision is
  $a_{11}^7a_{00}a_{30}a_{01}a_{12}^2-a_{11}^8a_{00}a_{30}a_{12}a_{02}=
  a_{11}^7a_{00}a_{30}a_{12}\cdot (a_{01}a_{12}-a_{11}a_{02})$. 
  The leading term of the coefficient for $f(x+t^b,y)$ of $y$ is $c_{11}
  t^b$, the leading term for $xy^2$ is $c_{12}t$, the one for $xy$ is
  $c_{11}$, and the one for $y^2$ is $c_{12}t^{1+b}$. Plugging those
  into $a_{01}a_{12}-a_{11}a_{02}$ yields $c_{11} t^b\cdot c_{12}t-
  c_{11}\cdot c_{12}t^{1+b}=0$.  
  Thus the generic valuation of $j$ (whose negative is equal to the
  decreasing cycle length) is not equal to the valuation of the
  $j$-invariant for values of $b$ in a whole interval. 

  We computed this example (as well as many other examples) using the
  procedure \texttt{drawtropicalcurve} from the 
  \textsc{Singular} library \texttt{tropical.lib} (see \cite{JMM07a}) which
  can be obtained via the URL
  \begin{center}
    http://www.mathematik.uni-kl.de/\textasciitilde keilen/en/tropical.html.
  \end{center}
  This library contains also a procedure \texttt{tropicalJInvariant}
  which computes the cycle length of a tropical curve as defined in
  Definition \ref{def-cyclelength}.
\end{remark}


\section{$\Delta$- equivalent marked subdivisions}\label{sec-delta}

In this section we want to show that the function ``cycle length'', $\cyclelength$, is
linear on the union of cones of the secondary fan of $\ca_3$ which
are $\Delta$-equivalent. Also, we provide the
classification of the different cases we need to consider in order to
compare the two linear functions $\val_\cdot(j)$ and $\cyclelength$ on 
such a union.  This is part of our first proof of Theorem
\ref{thm-main}.

\begin{remark}\label{rem-deltaequivalence}
  The Prime Factorization Theorem, \cite[Chap.\ 10, Thm.\ 1.2]{GKZ}, or
  alternatively  Lemma \ref{discriminant} tells us that the codimension $1$ cones of
  the Gr\"obner fan of $\Delta$ do not meet the interior of any
  top-dimensional cone of the secondary fan of $\ca_3$. Thus the Gr\"obner fan of
  $\Delta$ is a coarsening of the secondary fan of $\ca_3$. 
  Two cones of the secondary fan of $\ca_3$ are called
  \emph{$\Delta$-equivalent} if they are contained in the same cone of the Gr\"obner
  fan of $\Delta$. 

  It has been studied how two top-dimensional marked subdivisions
  whose cones belong to the same $\Delta$-equivalence class can differ. By
  \cite[Chap.\ 11, Prop.\ 3.8]{GKZ} they can be obtained from
  each other by a sequence of modifications along a circuit (see 
  \cite[Chap.\ 7, Sect. 2C]{GKZ}) such that each intermediate
  (top-dimensional) marked subdivision belongs to the same equivalence
  class. 
  Since our point configuration (i.e.\ $\ca_3$, the integer points of
  the triangle $Q_3$) is in the plane, we can use 
  \cite[Chap.\ 11, Prop.\ 3.9]{GKZ} to see that if a marked subdivision can be obtained from
  another equivalent one by a modification along a circuit, then this
  circuit consists of three collinear points on the boundary of
  $Q_3$. An example is shown in the following picture, the three points
  are $(0,0)$, $(1,0)$ and $(3,0)$. 

  \begin{center}
    \input{Graphics/deltaeq.pstex_t}
  \end{center}
\end{remark}

\begin{lemma}\label{lem-cyclinondelta}
  The function $\cyclelength$ (see Definition \ref{def-cyclelength}) is
  linear on a union of cones of the secondary fan of $\ca_3$ which
  are $\Delta$-equivalent (i.e.\ on a cone of the Gr\"obner
  fan of $\Delta$). 
\end{lemma}
\begin{proof}
  Given two $\Delta$-equivalent marked subdivisions $T$ and $T'$ of the
  secondary fan of $\ca_3$, we can use Lemma
  \ref{lem-cyclelengthlinear} to determine the function $\cyclelength$
  on the cone corresponding to each of them. Recall from Remark 
  \ref{rem-cyclelengthlinear} that the function is linear on each 
  cone of the secondary fan of $\ca_3$. We want to show that the
  assignment rules of these two linear
  functions coincide.
  
  Without restriction we can assume that $T$ can be obtained from $T'$
  by a modification along a circuit, and this circuit consists then of three collinear points on a
  facet of $Q_3$ (see Remark \ref{rem-deltaequivalence}).

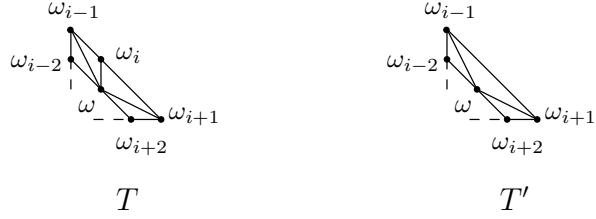
  \begin{center}
    \input{Graphics/circuit.pstex_t}
  \end{center}

  Recall from Lemma \ref{lem-cyclelengthlinear} that the coefficients of
  the linear function $\cyclelength$ can be determined using the
  determinants $D_{i,j}=\det(w_i,w_j)$, where $w_i=\omega_i-\tilde{\omega}$. 
  One easily sees that for $T$ and $T'$ the following two equations hold:
  \begin{eqnarray}\label{eq1}
    &D_{i-1,i}+D_{i,i+1}=D_{i-1,i+1}, \mbox{ and} \\[0.2cm]
    \label{eq2}
    & D_{i,i+1}=\lambda\cdot D_{i-1,i} \mbox{ for } \lambda \mbox{
      satisfying } \lambda\cdot (w_{i-1}-w_i)=w_{i}-w_{i+1}. 
  \end{eqnarray}

  To show that the two assignment rules of $\cyclelength$ on the cones
  for $T$ respectively $T'$ coincide we 
  have to show that for $T$ the summand for $\omega_i$ equals $0$ and
  the summand for $\omega_{i-1}$ equals the summand for $\omega_{i-1}$
  for $T'$. 
  The first statement follows immediately from Equation \eqref{eq1} above.
  To show the second statement, we subtract the two summands from each other:
  \begin{displaymath}
    \frac{D_{i-2,i-1}+D_{i-1,i}+D_{i,i-2}}{D_{i-2,i-1}\cdot
      D_{i-1,i}}-\frac{D_{i-2,i-1}+D_{i-1,i+1}+D_{i+1,i-2}}{D_{i-2,i-1}\cdot D_{i-1,i+1}} 
  \end{displaymath}

  Multiplied with $(D_{i-1,i}\cdot D_{i-1,i+1})$ this difference is equal to:
  \begin{multline*}
    D_{i-2,i-1}\cdot D_{i-1,i+1}+ D_{i-1,i}\cdot
    D_{i-1,i+1}+D_{i,i-2}\cdot D_{i-1,i+1}\\
    -D_{i-2,i-1}\cdot D_{i-1,i}-D_{i-1,i+1}\cdot
    D_{i-1,i}-D_{i+1,i-2}\cdot D_{i-1,i} \\ 
    = D_{i-2,i-1}\cdot D_{i,i+1}+D_{i,i-2}\cdot D_{i,i+1}+D_{i,i-2}\cdot
    D_{i-1,i} -D_{i+1,i-2}\cdot D_{i-1,i}\\ 
    = -\det(w_{i-1}-w_i,w_{i-2})\cdot
    D_{i,i+1}+\det(w_i-w_{i+1},w_{i-2})\cdot D_{i-1,i}  
    =0
  \end{multline*}
  where the 
  first equality follows from Equation \eqref{eq1} above and the last from \eqref{eq2}. 
\end{proof}

\begin{definition}
  Let us fix a cone $C_T$ of the secondary fan of $\ca_3$ corresponding to the
  marked subdivision $T$. We then
  denote by $\eta_T(i,j)$ the coefficient of $u_{ij}$ in the assignment
  rule of the linear function $u\mapsto \val_u(\Delta)$ on $C_T$, 
  and by $c_T(i,j)$ we denote the coefficient of $u_{ij}$ in the
  assignment rule of the linear function $\cyclelength$ restricted to
  $C_T$.
\end{definition}

\begin{remark}\label{rem-eta}  
  Note that by Lemma \ref{lem-genvallinear} and Remark \ref{rem-deltaequivalence}
  $\eta_T(i,j)=\eta_{T'}(i,j)$ for all $(i,j)\in\ca_3$ whenever $T$ and
  $T'$ belong to $\Delta$-equivalent cones of the secondary fan of
  $\ca_3$, and by Lemma \ref{lem-cyclinondelta} also $c_T(i,j)=c_{T'}(i,j)$ for all
  $(i,j)\in\ca_3$ in this situation.
\end{remark}

\begin{lemma}\label{lem-comparelin}
  Let $T$ be a marked subdivision of $(Q_3,\ca_3)$ corresponding to a
  top-dimensional cone in the secondary fan of $\ca_3$ (i.e.\ a
  triangulation) such that $(1,1)$ is a vertex of some marked polytope
  in $T$ (i.e.\ all dual plane tropical curves have a cycle). 
  Then $c_T(1,1)=\eta_T(1,1)-12$ and $c_T(i,j)=\eta_T(i,j)$ for all $(i,j)\neq (1,1)$. 
\end{lemma}
\begin{proof}
  Due to Remark \ref{rem-eta} we may for the proof assume that
  $T=\{(Q_\theta,\ca_\theta)\;|\;\theta\in\Theta\}$ is 
  the representative of its 
  $\Delta$-equivalence class with as few edges as possible.

  Moreover, if two triangulations $T$ and $T'$ can be transformed
  into each other by an integral unimodular linear isomorphism, i.e.\ by
  linear coordinate change of the projective coordinates $(x,y,z)$ with a
  matrix in $\Gl_3(\Z)$, and the claim holds for $T$ then it obviously
  also holds for $T'$. In this situation we say that $T$ and $T'$ are
  \emph{symmetric} to each other. We therefore only have to prove the
  claim \emph{up to symmetry}. 

  We want to use  \cite[Chap.\ 11, Thm.\ 3.2]{GKZ} which explains how 
  $\eta_T(i,j)$ can be computed. For each $(i,j)\in\ca_3$
  we have to consider all $(Q_\theta,\ca_\theta)$ such that
  $(i,j)\in\ca_\theta$. Note that since $T$ by assumption is a
  triangulation then $(i,j)\in\ca_\theta$ implies necessarily that
  $(i,j)$ is a vertex of $Q_\theta$. We have to distinguish four cases, where in
  the formulas $\vol(Q_\theta)$ denotes he generalized lattice volume
  (i.e.\ twice the euclidean area of $Q_\theta$):
  \begin{itemize}
  \item If $(i,j)$ is a vertex of $Q_3$, then 
    $\eta_T(i,j)=1-l_1-l_2+\sum_{(i,j)\in \ca_\theta} \vol(Q_\theta)$ where
    $l_1$ and $l_2$ denote the lattice lengths of those facets of some
    $Q_\theta$ adjacent to $(i,j)$ which are contained in facets of
    $Q_3$. E.g.\ if $(i,j)=(0,3)$ in the following triangulation $T$,
    then
    $\eta_T(0,3)=1-l_1-l_2+\vol(Q_{\theta_1})+\vol(Q_{\theta_2})=1-3-2+3+2=1$. 
    
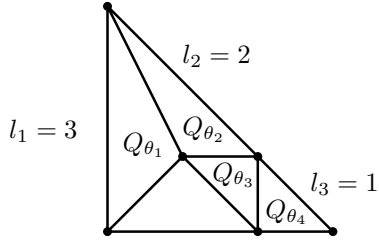
\begin{center}
      \begin{texdraw}
        \drawdim cm  \relunitscale 1 \linewd 0.03 
        \setgray 0.9
        \lfill f:0.9
        \setgray 0
        \move (0 0) \lvec (3 0) \lvec(0 3) \lvec (0 0) \lvec (1 1)
        \lvec (0 3)
        \move (2 0) \lvec (2 1) \lvec (1 1) \lvec (2 0)
        \move (0 0) \fcir f:0 r:0.06
        \move (0 3) \fcir f:0 r:0.06
        \move (3 0) \fcir f:0 r:0.06
        \move (1 1) \fcir f:0 r:0.06
        \move (2 0) \fcir f:0 r:0.06
        \move (2 1) \fcir f:0 r:0.06
        \htext (0.2 1){$Q_{\theta_1}$}
        \htext (1 1.2){$Q_{\theta_2}$}
        \htext (1.4 0.6){$Q_{\theta_3}$}
        \htext (2.1 0.1){$Q_{\theta_4}$}
        \htext (1 2.2){$l_2=2$}
        \htext (-1.3 1.2){$l_1=3$}
        \htext (2.7 0.5){$l_3=1$}
      \end{texdraw}
    \end{center}
  \item If $(i,j)$ lies on a facet of $Q_3$, is not a vertex of
    $Q_3$, but is a vertex of some $Q_{\theta'}$, then
    $\eta_T(i,j)=-l_1-l_2+\sum_{(i,j)\in \ca_\theta} \vol(Q_\theta)$ where
    again $l_1$ and $l_2$ denote the lattice lengths of those facets of some
    $Q_\theta$ adjacent to $(i,j)$ which are contained in facets of
    $Q_3$, e.g.\ if in the previous example $(i,j)=(2,1)$ then
    $\eta_T(i,j)=-l_2-l_3+\vol(Q_{\theta_2})+\vol(Q_{\theta_3})+\vol(Q_{\theta_4})=-2-1+2+1+1=1$. 
  \item If $(i,j)$ lies on a facet of $Q_3$, is not a vertex of
    any $Q_\theta$, then $\eta_T(i,j)=0$.
  \item And finally $\eta_T(1,1)=\sum_{(1,1)\in \ca_\theta} \vol(Q_\theta)$.
  \end{itemize}

  Let $Q$ be the union of all those $Q_\theta$ which contain $(1,1)$,
  and endow the marked polytope $(Q,Q\cap\ca_3)$ 
  with the subdivision, say $T_Q$, induced by $T$.
  We say that $Q$ meets a facet of $Q_3$ if the intersection of $Q$ with
  this facet is $1$-dimensional (and not only a vertex). Moreover, we say that a
  facet of $Q$ is multiple if it contains more than two lattice
  points. 

  We first want to show that $\eta_T(i,j)$ and $c_T(i,j)$ are as claimed
  whenever $(i,j)\in Q$.
  Up to symmetry, we have to distinguish the following cases for $Q$ and
  $T_Q$:
  \begin{itemize}
  \item Assume $Q$ meets all three facets of $Q_3$, and assume that for
    all three facets  the intersection with $Q$ is multiple. Then $Q$
    looks (up to symmetry) like one of the following two pictures: 
    \begin{center}
      \input{Graphics/L1.pstex_t}
    \end{center}
    In the second case, $\eta_T(1,1)=8$. Using Lemma
    \ref{lem-cyclelengthlinear} we can compute $c_T(1,1)$. It is a sum
    with a summand for each vertex of $Q$. The summand for $(0,0)$ is  
    \begin{displaymath}
      \frac{\det\begin{pmatrix}1&-1\\
          -1&-1\end{pmatrix}+\det\begin{pmatrix}-1&-1\\
          -1&2\end{pmatrix}+\det\begin{pmatrix}-1&1\\
          2&-1\end{pmatrix}}{\det\begin{pmatrix}1&-1\\
          -1&-1\end{pmatrix}\cdot\det\begin{pmatrix}-1&-1\\
          -1&2\end{pmatrix}}=-1.
    \end{displaymath}
    Computing the other $3$ summands analogously we get $c_T(1,1)=-4=\eta_T(1,1)-12$.
    In the first case, $\eta_T(1,1)=9$ and $c_T(1,1)=-3$. 
  \item Assume $Q$ meets two facets of $Q_3$ multiply and one facet non-multiply.
    \begin{center}
      \input{Graphics/L2.pstex_t}
    \end{center}
    In both cases, $\eta_T(1,1)=7$ and $c_T(1,1)=-5$. 
  \item Assume $Q$ meets two facets of $Q_3$ multiply and the third
    facet not at all.
    \begin{center}
      \input{Graphics/L3.pstex_t}
    \end{center}
    In both cases, $\eta_T(1,1)=6$ and $c_T(1,1)=-6$. 
  \end{itemize}
  \begin{itemize}
  \item Assume $Q$ meets only one facet of $Q_3$ multiply (and the two
    remaining facets non-multiply, or only one of them and that one non-multiply, or none of them
    at all).\begin{center}
      \input{Graphics/L4.pstex_t}
    \end{center}
    In the first case, $\eta_T(1,1)=6$ and $c_T(1,1)=-6$, in the second and
    third case, $\eta_T(1,1)=5$ and $c_T(1,1)=-7$, and in the last case,
    $\eta_T(1,1)=4$ and $c_T(1,1)=-8$.  
  \item Assume $Q$ meets 3 facets of $Q_3$, but none of them multiply.
    \begin{center}
      \input{Graphics/L5.pstex_t}
    \end{center}
    In the first case, $\eta_T(1,1)=5$ and $c_T(1,1)=-7$, and in the second 
    case, $\eta_T(1,1)=6$ and $c_T(1,1)=-6$.
  \item Assume $Q$ meets only two facets of $Q_3$, and none of them multiply.
    \begin{center}
      \input{Graphics/L6.pstex_t}
    \end{center}
    In the first case, $\eta_T(1,1)=5$ and $c_T(1,1)=-7$, and in the second 
    case, $\eta_T(1,1)=4$ and $c_T(1,1)=-8$.
  \item Assume $Q$ meets only one facet of $Q_3$ and it does so non-multiply.
    \begin{center}
      \input{Graphics/L7.pstex_t}
    \end{center}
    In the first and second case, $\eta_T(1,1)=4$ and $c_T(1,1)=-8$, and in the third 
    case, $\eta_T(1,1)=3$ and $c_T(1,1)=-9$.
  \item Assume $Q$ meets no facet of $Q_3$ at all.
    \begin{center}
      \input{Graphics/L8.pstex_t}
    \end{center}
    Finally, in this case, $\eta_T(1,1)=3$ and $c_T(1,1)=-9$.
  \end{itemize}
  Thus the claim for $(1,1)$ is shown.
  Now assume $(1,1)\not=(i,j)\in Q$ is not a vertex of $Q_3$. 
  If $(i,j)$ is also not a vertex of any $Q_\theta$ then there is no
  edge in the subdivision from $(1,1)$ to $(i,j)$ and thus $(i,j)$ does
  not contribute to the formula for the cycle length, i.e.\
  $c_T(i,j)=0$. However, the same holds true also for $\eta_T(i,j)$. We
  may thus assume that $(i,j)$ is a vertex of some $Q_\theta$, and we
  may without restriction assume $(i,j)=(0,1)$. The
  classification of cases we have to consider is very similar to the
  above, and we will not give the details -- in particular, leaving the
  computation of $c_T(i,j)$ and $\eta_T(i,j)$ to the reader. We do not
  have to consider 
  the whole of $Q$, but only the triangles which are adjacent to $(i,j)$. 
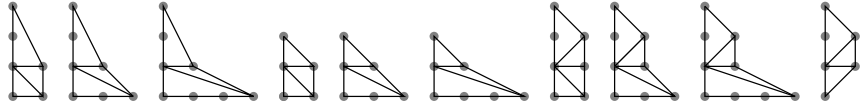
  \begin{center}
    \input{Graphics/boundary.pstex_t}
  \end{center}
  If $(1,1)\not=(i,j)\in Q$ is a vertex  of $Q_3$ (without restriction
  $(i,j)=(0,0)$), the following cases have to be considered: 
  \begin{center}
    \input{Graphics/vertex.pstex_t}
  \end{center}
  Finally, we have to consider the case were $(i,j)$ is not part of
  $Q$. Obviously, $c_T(i,j)=0$ in this case and we have to show the same
  for $\eta_T(i,j)$. 
  Assume first that $(i,j)$ is a vertex of $Q_3$, without restriction we can assume $(i,j)=(0,0)$.
  There must be a facet of $Q$ such that $(0,0)$ is on one side of it and $(1,1)$ is on the other.
  Then (up to symmetry) there are $3$ possibilities for that facet.
  \begin{center}
    \input{Graphics/Lvert.pstex_t}
  \end{center}
  Since we assumed that $T$ is the representative with as few edges as
  possible, the triangle formed by that facet of $Q$ and $(0,0)$ can not be
  additionally subdivided in the second and third picture. In any case,
  $(0,0)$ is a vertex of only one triangle, which has one facet of
  integer length $1$ and one facet of
  integer length $l$ where $1\leq l\leq3$. Thus
  $\eta_T(0,0)=1-1-l+l=0$. 
  Now assume that $(i,j)$ is not a vertex of $Q_3$, without restriction
  $(i,j)=(1,0)$. Again there must be a facet of $Q$ such that $(1,0)$ is on one
  side and $(1,1)$ is on the other. Up to symmetry this can only be one of the line segments
  in the two right pictures above. We assumed that $T$ is the
  representative of its $\Delta$-equivalence class with as few edges as
  possible. But that means there is no edge through $(1,0)$ and $(1,0)$
  is not a vertex of a triangle in the subdivision. Thus
  $\eta_T(1,0)=0$. 
\end{proof}

\begin{remark}
  In the proof above the computation that shows that
  $\eta_T(i,j)=c_T(i,j)$ is different in each of the considered
  cases. In particular, in the computation of $\eta_T(i,j)$ the part of
  $Q_3$ which is not part of the cycle is involved while this is not the
  case for $c_T(i,j)$.
  Therefore it is most unfortunately  not possible to replace the
  consideration of several cases by an argument which holds for all of
  them at the same time.

  However, using \texttt{polymake} and \textsc{Singular} one can compute
  the vertices of the Newton polytope of $\Delta$ and for each vertex
  one can compute the dual cone in the Gr\"obner fan of $\Delta$ and the
  triangulation of $(Q_3,\ca_3)$ with as few edges as possible
  corresponding to this cone. That way one can verify the above
  computations for $c_T$ and $\eta_T$, since the values for $\eta_T$ can
  be read off immediately from the exponents of the vertex of the Newton
  polytope, while the $c_T$ can be computed with the formula in Lemma
  \ref{lem-cyclelengthlinear}. These computations have been made using
  the procedure \texttt{displayFan} and the 
  result can be obtained via the URL 
  \begin{center}
    http://www.mathematik.uni-kl.de/\textasciitilde keilen/en/jinvariant.html.
  \end{center}
  The advantage is that the file \texttt{discriminant\_fan\_of\_cubic.ps} 
  available via this url shows the cases not
  only up to symmetry, but it shows actually all possible cases.
\end{remark}

\section{Discriminant}\label{sec-discriminant}

The aim of this section is to show that the Gr\"obner fan of $\Delta$
is a coarsening of the secondary fan. 
This follows from the
Prime Factorization Theorem of Gelfand, Kapranov and Zelevinsky
(\cite[Chap.\ 10, Sect.\ 2]{GKZ}, see also  \cite[Conj.\ 5.2]{DFS05}).
We present our own proof because we hope that the techniques will be
useful.  They are similar to 
those of Sturmfels (see \cite{Stu94}).

Our reference for toric varieties and the $\ca$-discriminant is
\cite{GKZ}, and we will use Notation \ref{not-x}. Let $\ca\subset
\Z^\Lambda$ be a finite set of lattice points. It 
defines a projective toric variety $X_\ca\subset \P^{|\ca|-1}$ over
$K$, where $K$ is any field with non-archemedian valuation
$\val:K^*\rightarrow \R$ whose value group is dense in $\R$ with
respect to the Euclidean topology. 

Our special case is $\ca_3=Q_3\cap \Z^2$. In this case, $X_\ca=\P^2$
embedded in $\P^9$ by the $3$-uple embedding.   

We review some facts about toric varieties. A Laurent polynomial with support $\ca$,
\[f=\sum_{\omega\in\ca} a_\omega \underline{x}^\omega\]
can be thought of as defining a hypersurface in the toric variety
$X_\ca$ (see \cite{GKZ}), and since this hypersurface 
coincides with the hyperplane section defined by the coefficients
of $f$ under the embedding $X_\ca\subseteq\P^{|\ca|-1}$ we 
identify a polynomial $f=\sum_{\omega\in\ca}a_\omega\cdot
\underline{x}^\omega$ with the point $[a_\omega]_{\omega\in\ca}$ in
the dual space $(\P^{|\ca|-1})^\vee$ of $\P^{|\ca|-1}$.  

The \emph{dual variety} $X_\ca^\vee\subset (\P^{|\ca|-1})^\vee$ is the
Zariski closure of the set of all hyperplanes $H$ such that $X_\ca^\circ\cap H$ is singular where $X_\ca^\circ$ is the open torus in $X_\ca$.  If $X_\ca^\vee$
is of codimension $1$, then its defining polynomial is called the
\emph{discriminant} $\Delta_\ca$, otherwise we say that the $X_\ca$
has a \emph{degenerate dual variety} and we set $\Delta_\ca=1$.
For our special case $\ca_3$, $\Delta=\Delta_{\ca_3}$ 
is the denominator of the $j$-invariant.

\begin{lemma}\label{lem-smoothdelta}
  Suppose $X_\ca$ is smooth and $f\in V(\Delta_\ca)$, then the variety 
  $V(f)$, considered as a hypersurface in $X_\ca$, is singular. 
\end{lemma}

\begin{proof}
  Consider the universal hypersurface
  $\cu=\left\{\sum_{\omega\in\ca}a_\omega\cdot\underline{x}^\omega=0\right\}\subset
  (\P^{|\ca|-1})^\vee\times X_\ca$ and the flat subfamily defined by the
  vanishing of $\Delta_\ca$:
  \begin{displaymath}
    \xymatrix{
      \left\{(f,x)\in\cu\;|\;\Delta_\ca(f)=0\right\}
      \ar@{^{(}->}[r]\ar[dr] &V(\Delta_\ca)\times X_\ca\ar[d]\\
      &V(\Delta_\ca)
    }
  \end{displaymath}
  By definition the general fiber of this family is singular, and
  thus so must be each fiber.   
\end{proof}

The fact, however, that this singular point need not be in the open
torus is a problem when proving that the Gr\"obner fan of the
discriminant $\Delta_\ca$ is a coarsening of the secondary fan of
$\ca$. For this we will have to reduce to the restriction of $f$ to
some face $\Gamma$ which has a singular point in the torus orbit
$X_\Gamma^0$ (see Lemma \ref{lem-singularorbit}). 

Let ${Q_\ca}=\Conv(\ca)$ and suppose
$\ca={Q_\ca}\cap\Z^\Lambda$.   For each face $\Gamma$ of ${Q_\ca}$ of dimension
$k$, we have a parameterization of the open torus orbit  $X^0_\Gamma$
given by 
\[i_\Gamma:(K^*)^k\hookrightarrow X_\ca.\]
We may consider the restriction 
\[f_\Gamma=\sum_{\omega\in\ca\cap\Gamma} a_\omega \underline{x}^\omega\]
of $f$ to $X^0_\Gamma$ as a function 
\begin{displaymath}
  f_\Gamma:(K^*)^k\rightarrow K:\xi\mapsto f_\Gamma\big(i(\xi)\big).
\end{displaymath}
Picking coordinates $y_1,\dots,y_k$ on $(K^*)^k$, we define
the subset $Z_\Gamma$ of $(\P^{|\ca|-1})^\vee\times X_\ca$ by
\[Z_\Gamma=\bigcup_{f\in(\P^{|\ca|-1})^\vee}\{f\}\times
i_\Gamma\left(V(f_\Gamma)\cap V\left(\frac{\partial 
    f_{\Gamma}}{\partial y_1}\right)\cap\dots\cap
V\left(\frac{\partial f_{\Gamma}}{\partial y_k}\right)\right).\]
$Z_\Gamma$ can be viewed as the subset of $(\P^{|\ca|-1})^\vee\times X_\ca$
consisting of pairs $(f,x)$ of functions $f$ and points $x$ such that
$f_\Gamma$ is singular on $X^0_\Gamma$ at $x$.   Note that if $\Gamma$
is just a point, then $X^0_\Gamma$ is a closed point of $X_\ca$, and
$f$ is singular on $X^0_\Gamma$ if and only if it is zero on
$X^0_\Gamma$. 

\begin{lemma} \label{singularclosure}
  For two faces $\Gamma'$ and $\Gamma$ of $Q_\ca$ with
  $\Gamma'\subset\Gamma$ there is the following inclusion
  \begin{displaymath}
    \overline{Z_\Gamma}\cap \left((\P^{|\ca|-1})^\vee\times
      X_{\Gamma'}^0\right)
    \subseteq Z_{\Gamma'}^0,
  \end{displaymath}
  where $\overline{Z_\Gamma}$ denotes the Zariski closure of
  $Z_\Gamma$ in $(\P^{|\ca|-1})^\vee\times X_\ca$.
\end{lemma}

\begin{proof}
  We may restrict to the toric variety $X_\Gamma$ and therefore, we
  may suppose that $\Gamma=Q_\ca$.  By applying a monomial
  change-of-variables (i.e. a change of variables of the form
  $y_\nu=\prod_{\lambda\in\Lambda} x_\lambda^{a_{\nu\lambda}}$ for $\nu\in\Lambda$
  where $(a_{\nu\lambda})_{\nu,\lambda\in\Lambda}$
  is an integer matrix of determinant $\pm 1$), we may suppose that
  $\Gamma'$ lies in a 
  coordinate subspace of $\R^\Lambda$ given by
  $x_\lambda=0$ for $\lambda\in\Lambda\setminus\Lambda'$ with
  $\#\Lambda'=\dim(\Gamma')$. 

  By construction, there is an embedding
  $i:X_\ca\hookrightarrow\P^{|\ca|-1}$ of $X_\ca$ into projective
  space.  Let us pick homogeneous 
  coordinates $[b_\omega]_{\omega\in\ca}$ so that on the open torus
  $X_\ca^0=(K^*)^\Lambda$ the map $i$ is given by 
  \begin{displaymath}
    (K^*)^\Lambda\longrightarrow \P^{|\ca|-1}:\underline{x}\mapsto b_\omega=\underline{x}^\omega.
  \end{displaymath}
  Now consider the following morphism
  \begin{displaymath}
    j=\id\times i:(\P^{|\ca|-1})^\vee\times X_\ca\rightarrow
    (\P^{|\ca|-1})^\vee\times \P^{|\ca|-1},
  \end{displaymath}
  which is the identity on the first factor and $i$
  on the second factor.  We denote the coordinates on 
  $(\P^{|\ca|-1})^\vee$ again by $[a_\omega]_{\omega\in\ca}$.  There is a
  universal hypersurface 
  $\cv\subset (\P^{|\ca|-1})^\vee\times \P^{|\ca|-1}$ cut out by 
  \begin{displaymath}
    \sum_{\omega\in\ca} a_\omega \cdot b_\omega=0.
  \end{displaymath}
  Observe that $\cu=j^{-1}(\cv)$ is the universal
  hypersurface on $(\P^{|\ca|-1})^\vee\times X_\ca$ defined by the vanishing
  of $\sum_{\omega\in\ca}a_\omega\cdot\underline{x}^\omega$.  We may also
  define a universal singular locus $\cv_{\Sing}\subset
  (\P^{|\ca|-1})^\vee\times\P^{|\ca|-1}$ by the $\#\Lambda+1$ equations 
  \begin{displaymath}
    \sum_{\omega\in\ca} a_\omega\cdot b_\omega=
    \sum_{\omega\in\ca} \omega_\lambda\cdot a_\omega\cdot b_\omega=0
    \;\;\;\mbox{ for }\;  \lambda\in\Lambda.
  \end{displaymath} 
  It is straight forward from the definitions that 
  \begin{displaymath}
    j^{-1}(\cv_{\Sing})\cap
    \big((\P^{|\ca|-1})^\vee\times X_\Gamma^0\big)=Z_\Gamma.
  \end{displaymath}   
  Therefore, $\ol{Z_\Gamma}\subseteq j^{-1}(\cv_{\Sing})$.

  We may define 
  $\cv_{\Gamma',\Sing}\subseteq
  (\P^{|\ca|-1})^\vee\times\P^{|\ca|-1}$ by 
  \begin{displaymath}
    \sum_{\omega\in\Gamma'\cap\ca} a_\omega\cdot b_\omega=
    \sum_{\omega\in\Gamma'\cap\ca} \omega_\lambda\cdot a_\omega\cdot b_\omega
    =0 
    \;\;\;\mbox{ for }\;\lambda\in\Lambda'.
  \end{displaymath}
  Again we get immediately from the definition that
  \begin{displaymath}
    Z_{\Gamma'}=
    j^{-1}(\cv_{\Gamma',\Sing})\cap \big((\P^{|\ca|-1})^\vee\times X_{\Gamma'}^0\big)=
    j^{-1}(\cx)\cap \big((\P^{|\ca|-1})^\vee\times X_{\Gamma'}^0\big),
  \end{displaymath}
  where
  \begin{displaymath}
    \cx=\cv_{\Gamma',\Sing}\cap\big(\big(\P^{\ca|-1}\big)^\vee\times 
    V(b_\omega\;|\;\omega\not\in\Gamma')\big)
    =\cv_{\Sing}\cap\big((\P^{|\ca|-1})^\vee\times V(b_\omega\;|\;\omega\not\in\Gamma')\big).
  \end{displaymath}
  Thus by taking inverse images by $j$, we get 
  \begin{multline*}
    \ol{Z_\Gamma}\cap\big((\P^{|\ca|-1})^\vee\times X_{\Gamma'}^0\big)\subseteq
    j^{-1}(\cv_{\Sing}) \cap\big((\P^{|\ca|-1})^\vee\times X_{\Gamma'}^0\big)\\
    =j^{-1}(\cx) \cap\big((\P^{|\ca|-1})^\vee\times X_{\Gamma'}^0\big)
    =Z_{\Gamma'}.    
  \end{multline*}   
\end{proof}

\begin{lemma} \label{lem-singularorbit}
  If $f\in V(\Delta_\ca)$, then
  there is some face $\Gamma$ of $Q_\ca$ so that 
  $f_\Gamma$ is singular on $X^0_\Gamma$. 
\end{lemma}

\begin{proof}  
  If $\cf({Q_\ca})$ denotes the set of all faces of ${Q_\ca}$ including ${Q_\ca}$
  itself, we have to show that 
  \begin{equation}\label{eq:singularorbit:1}
    V(\Delta_\ca)\subseteq \bigcup_{\Gamma\in\cf({Q_\ca})}\pi(Z_\Gamma)
    =\pi\left(\bigcup\nolimits_{\Gamma\in\cf({Q_\ca})} Z_\Gamma\right),
  \end{equation}
  where $\pi:(\P^{|\ca|-1})^\vee\times X_\ca\longrightarrow(\P^{|\ca|-1})^\vee$
  denotes the projection onto the first factor. Since by definition
  $V(\Delta_\ca)$ is the Zariski closure of $\pi(Z_{Q_\ca})$ it suffices to
  show that the right hand side in \eqref{eq:singularorbit:1} is
  Zariski closed, or equivalently that $\bigcup_{\Gamma\in\cf({Q_\ca})}
  Z_\Gamma$ is so. This, however, follows once we know that 
  \begin{equation}
    \label{eq:singularorbit:2}
    \overline{Z_\Gamma}\subseteq \bigcup_{\Gamma'\in\cf({Q_\ca})}Z_{\Gamma'}
  \end{equation}
  for all $\Gamma\in\cf({Q_\ca})$, where $\overline{Z_\Gamma}$ denotes the Zariski closure of
  $Z_\Gamma$.   Since the Zariski closure of $X_\Gamma^0$ is 
  $\overline{X_\Gamma^0}=\bigcup_{\Gamma'\in\cf(\Gamma)}
  X_{\Gamma'}^0$ and $\overline{Z_\Gamma}\subseteq
  \overline{X_\Gamma^0}$, this fact follows from the Lemma \ref{singularclosure}.
\end{proof}

\begin{lemma}\label{lem-singularsimplex}
  Let $g\in\C[y_1,\ldots,y_k]$ be a polynomial with $k+1$ terms whose 
  Newton polytope $N(g)$ is a $k$-dimensional simplex.  
  Then the hypersurface
  $V(g)$ has no singular points in $(\C^*)^k$.
\end{lemma}

\begin{proof}
  We may divide $g$ by a monomial without changing
  $V(g)\cap(\C^*)^k$.  Therefore, we may suppose that one vertex of
  $N(g)$ is at 
  the origin.  By applying a monomial change-of-variables, i.e.\ a
  coordinate change of the form $y_i\mapsto y_1^{a_{i1}}\cdots
  y_k^{a_{ik}}$ with $a_{ij}\in\Z$ and
  $\det\big((a_{ij})_{i,j=1,\ldots,k}\big)=\pm 1$, we may 
  suppose that the edges from $0$ to the other $k$ vertices of $N(g)$ lie
  along the axes.  Therefore,
  \[g=a+\sum_{i=1}^k a_i y_i^{c_i}\]
  for $c_i\in\N$. In that case $\frac{\partial f}{\partial
    y_i}$ is a monomial and has no root in $(\C^*)^\Lambda$.
\end{proof}

\begin{proposition}\label{discriminant}
  Let $\ca\subset\Z^\Lambda$ be such that $X_\ca^\vee$ has codimension
  one and the discriminant $\Delta_\ca$
  exists. Then the tropicalization 
  $\Trop\big(V(\Delta_\ca)\big)$ of $V(\Delta_\ca)$ is supported on the codimension
  one cones of the secondary fan of $\ca$ in $\R^\ca$.
\end{proposition}
\begin{proof}
  We assume that $\Trop\big(V(\Delta_\ca)\big)$ intersects the
  interior of a top dimensional cone of the secondary fan of $\ca$ in
  a point $u\in\Q^\ca$ and derive a contradiction.  

  Since $u\in\Trop\big(V(\Delta_\ca)\big)$ by the Lifting Lemma for
  hypersurfaces (see e.g.\ \cite[Thm.~2.13]{EKL04}) we can lift $u$ to a point
  $f=[a_\omega]_{\omega\in\ca}\in V(\Delta_\ca)$ such that 
  $\val(a_\omega)=u_\omega$ for all $\omega\in\ca$ -- in particular,
  \begin{equation}
    \label{eq:aomega}
    a_\omega\not=0\;\;\; \mbox{ for all }\;\omega\in\ca.
  \end{equation} 
  By Lemma \ref{lem-singularorbit} there exists then a, say
  $k$-dimensional, face $\Gamma$ of ${Q_\ca}=\Conv(\ca)$ such that
  $V(f_\Gamma)\subset(K^*)^k$ has a singular point $\xi\in(K^*)^k$, and
  $k>0$ due to \eqref{eq:aomega}.
  We define $\omega=\big(\val(\xi_1),\ldots,\val(\xi_k)\big)$ to be
  the valuation of $\xi$, and we then claim that
  $\xi_0=\big(\lc(\xi_1),\ldots,\lc(\xi_k)\big)\in(\C^*)^k$ is a singular
  point of $\tini_\omega(f_\Gamma)$. 
  In order to see that $\xi_0$ is a singular point of
  $\tini_\omega(f_\Gamma)$ it suffices to note that 
  \begin{displaymath}
    \tini_\omega\left(\frac{\partial f_\Gamma}{\partial x_i}\right)
    =
    \frac{\partial \tini_\omega(f_\Gamma)}{\partial x_i},
  \end{displaymath}
  and that for any polynomial $g\in K[y_1,\ldots,y_k]$ with $g(\xi)=0$
  we necessarily have $\tini_\omega(g)(\xi_0)=0$.

  Since $u$ is in the interior of a full-dimensional cone of the
  secondary fan of $\ca$ the Newton subdivision, say
  $\{(Q_\theta,\ca_\theta)\;|\;\theta\in\Theta\}$, of $f$ is a
  triangulation. By definition the Newton polytope of the $t$-initial form
  $\tini_\omega(f_\Gamma)$ is a face of some $Q_\theta$ and is thus a
  simplex. But then by Lemma \ref{lem-singularsimplex}
  $\tini_\omega(f_\Gamma)$ has no singular point in the torus
  $(\C^*)^k$ in contradiction to the existence of $\xi_0$.
\end{proof}


\section{Numerator of the $j$-invariant}\label{sec-num}

Unfortunately, for the numerator $A$ of the $j$-invariant it is not
true that the Gr\"obner fan of $A$ is 
a coarsening of the secondary fan, as follows from Example
\ref{ex-notcoarse}.

\begin{example}\label{ex-notcoarse}
  We provide an example which shows that the Gr\"obner fan of $A$ is not
  a coarsening of the secondary fan in the case of curves of a
  particular 
  form.  The case of the full cubic is more complicated but
  analogous. It can easily be proved by a computation using
  \texttt{polymake} -- this can be done using the procedure \texttt{nonrefinementC} in the
  library \texttt{jinvariant.lib} (see \cite{KMM07a}). 

  Let us consider curves in Weierstrass form
  \[y^2+axy-x^3-bx^2-1=0.\]
  This corresponds to taking $\ca=\{(0,2),(1,1),(3,0),(2,0),(0,0)\}$.
  The fixing the constant coefficient and the coefficients of  $y^2$
  and $x^3$ has the effect of fixing an isomorphism
  $\R^2\cong R^{|\ca|}/L$ in light of Remark \ref{linearityspace}.   
    
  By the usual formulas for the $j$-invariant, we have
  \begin{displaymath}
    A=(a^2+4b)^6\;\;\;\mbox{ and }\;\;\; \Delta=-(a^2+4b)^3-432,
  \end{displaymath}
  so that
  \[j=-\frac{(a^2+4b)^6}{(a^2+4b)^3+432}.\]

  The following picture shows the tropicalization of the numerator
  $A$, the tropicalization of the denominator $\Delta$, and the
  secondary fan in $\R^{|\ca|}/L$. 

  
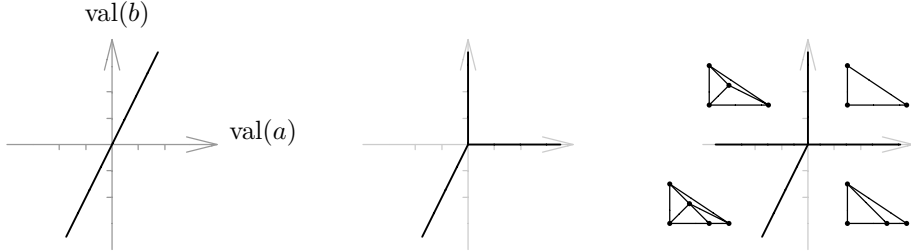
\bigskip
  \begin{center}
    \begin{texdraw}
      \drawdim cm \relunitscale 0.35 \linewd 0.05 \lpatt (1 0) \setgray 0.6 \arrowheadtype t:V
      \move (-4 0) \avec (4 0) \move (0 -4) \avec (0 4)
      \move (-2 0) \rlvec (0 -0.2)
      \move (-1 0) \rlvec (0 -0.2)
      \move (1 0) \rlvec (0 -0.2)
      \move (2 0) \rlvec (0 -0.2)
      \move (0 -2) \rlvec (-0.2 0)
      \move (0 -1) \rlvec (-0.2 0)
      \move (0 1) \rlvec (-0.2 0)
      \move (0 2) \rlvec (-0.2 0)
      \htext (4.5 0) {$\val(a)$}
      \htext (-1 4.5) {$\val(b)$}
      \setgray 0 \linewd 0.07
      \move (-1.75 -3.5) \lvec (1.75 3.5) 
    \end{texdraw}
    \hspace{1.5cm}
    \begin{texdraw}
      \drawdim cm \relunitscale 0.35 \linewd 0.05 \lpatt (1 0) \setgray 0.8 \arrowheadtype t:V
      \move (-4 0) \avec (4 0) \move (0 -4) \avec (0 4)
      \move (-2 0) \rlvec (0 -0.2)
      \move (-1 0) \rlvec (0 -0.2)
      \move (1 0) \rlvec (0 -0.2)
      \move (2 0) \rlvec (0 -0.2)
      \move (0 -2) \rlvec (-0.2 0)
      \move (0 -1) \rlvec (-0.2 0)
      \move (0 1) \rlvec (-0.2 0)
      \move (0 2) \rlvec (-0.2 0)
      \setgray 0 \linewd 0.07
      \move (-1.75 -3.5) \lvec (0 0) \lvec (0 3.5) 
      \move (0 0) \lvec (3.5 0)
    \end{texdraw}
    \hspace{1cm}
    \begin{texdraw}
      \drawdim cm \relunitscale 0.35 \linewd 0.05 \lpatt (1 0) \setgray 0.8 \arrowheadtype t:V
      \move (-4 0) \avec (4 0) \move (0 -4) \avec (0 4)
      \move (-2 0) \rlvec (0 -0.2)
      \move (-1 0) \rlvec (0 -0.2)
      \move (1 0) \rlvec (0 -0.2)
      \move (2 0) \rlvec (0 -0.2)
      \move (0 -2) \rlvec (-0.2 0)
      \move (0 -1) \rlvec (-0.2 0)
      \move (0 1) \rlvec (-0.2 0)
      \move (0 2) \rlvec (-0.2 0)
      \setgray 0 \linewd 0.07
      \move (-1.75 -3.5) \lvec (0 0) \lvec (0 3.5) 
      \move (-3.5 0) \lvec (3.5 0) \linewd 0.05      
      \relunitscale 1.5
      \move (-2.5 1) \rlvec (1.5 0) \rlvec (-1.5 1) \rlvec (0 -1) 
      \rlvec (0.5 0.5) \rlvec (1 -0.5) 
      \move (-2.5 2) \rlvec (0.5 -0.5)
      \move (-2.5 1) \fcir f:0 r:0.07
      \move (-1 1) \fcir f:0 r:0.07
      \move (-2 1.5) \fcir f:0 r:0.07
      \move (-2.5 2) \fcir f:0 r:0.07
      \move (-3.5 -2) \rlvec (1.5 0) \rlvec (-1.5 1) \rlvec (0 -1) 
      \rlvec (0.5 0.5) \rlvec (1 -0.5) 
      \move (-3.5 -1) \rlvec (0.5 -0.5) \rlvec (0.5 -0.5) 
      \move (-3.5 -2) \fcir f:0 r:0.07
      \move (-3.5 -1) \fcir f:0 r:0.07
      \move (-2 -2) \fcir f:0 r:0.07
      \move (-3 -1.5) \fcir f:0 r:0.07
      \move (-2.5 -2) \fcir f:0 r:0.07
      \move (1 1) \rlvec (1.5 0) \rlvec (-1.5 1) \rlvec (0 -1) 
      \move (1 1) \fcir f:0 r:0.07
      \move (2.5 1) \fcir f:0 r:0.07
      \move (1 2) \fcir f:0 r:0.07
      \move (1 -2) \rlvec (1.5 0) \rlvec (-1.5 1) \rlvec (0 -1) 
      \move (1 -1) \rlvec (1 -1) 
      \move (1 -2) \fcir f:0 r:0.07
      \move (1 -1) \fcir f:0 r:0.07
      \move (2.5 -2) \fcir f:0 r:0.07
      \move (2 -2) \fcir f:0 r:0.07
    \end{texdraw}
  \end{center}

  Observe that the tropicalization of the denominator is supported on
  the codimension one skeleton of the secondary fan while that of the
  numerator intersects a top-dimensional cone of the secondary fan. 
\end{example}

However, we are only interested in plane tropical cubics which have a cycle,
that is, which are dual to marked subdivisions for which the interior
point can be seen. 
All these cones of the secondary fan are completely contained in one
cone of the Gr\"obner fan of $A$. We verified this computationally
using \texttt{polymake} (see \cite{Pol97}). As usual we use the
coordinates $u_{ij}$ with $(i,j)\in\ca_3$ on $\R^{\ca_3}$ and we
denote by $e_{kl}=(\delta_{ik}\cdot\delta_{jl}\;|\;(i,j)\in\ca_3)$ the
canonical basis vector in $\R^{\ca_3}$ having a one in position $kl$
and zeros elsewhere.

\begin{lemma}\label{lem-num}
  Let $U$ be the union of all cones of the secondary fan of $\ca_3$ corresponding
  to marked subdivisions
  $T=\{(Q_\theta,\ca_\theta)\;|\;\theta\in\Theta\}$ of $(Q_3,\ca_3)$ for
  which $(1,1)$ is a vertex of some $Q_\theta$. Then $U$ 
  is contained in a single cone of the Gr\"obner fan of the $A$, namely in
  the cone dual to the vertex $12e_{11}$ of the Newton polytope of $A$. 
\end{lemma}

\begin{proof}
  As input for \texttt{polymake} we use all exponents of the polynomial
  $A\in\Q[\underline{a}]$. The convex hull of the set of all exponents
  is the Newton polytope, say $N(A)$,
  of $A$ and its vertices are the output of \texttt{polymake}. The
  Newton polytope has $19$ vertices. Dual to each vertex is a
  top-dimensional cone of the Gr\"obner fan $\mathcal{F}(A)$ of $A$, because the
  Gr\"obner fan is dual to the Newton polytope (see \cite[Thm.\ 2.5 and
  Prop.\ 2.8]{Stu96}). The inequalities describing the cone
  $C$ dual to the vertex $V$ are given by the hyperplanes orthogonal to
  the edge vectors connecting $V$ with its neighboring vertices in $N(A)$. We
  compute the neighboring vertices for each vertex using
  \texttt{polymake} and deduce thus inequalities for each of the
  top-dimensional cones of the Gr\"obner fan of $A$. 

  We do these computations identifying $\R^{\mathcal{A}_3}$ with
  $\R^{10}$ via the following ordering of the variables:
  \begin{displaymath}
    u_{11},u_{30},u_{20},u_{10},u_{00},u_{21},u_{01},u_{12},u_{02},u_{03}.
  \end{displaymath}
  In order for a marked subdivision $T_\psi=\{(Q_\theta,\mathcal{A}_\theta)\;\theta\in\Theta\}$ of
  $(Q_3,\mathcal{A}_3)$ given by 
  $$\psi:\R^{\mathcal{A}_3}\rightarrow \R:(i,j)\mapsto u_{ij}$$
  to have the point $(1,1)$ as vertex of some $Q_\theta$ it is obviously
  necessary that the $u_{ij}$ satisfy the following inequalities:
  \begin{equation}\label{eq-inequalities}
    \begin{array}{rclp{0.3cm}rcl}
      3\cdot u_{01} + 2\cdot u_{30} + u_{03} & > & 6 \cdot u_{11} &&
      3\cdot u_{10} + 2\cdot u_{03} + u_{30} & > & 6 \cdot u_{11}\\
      3\cdot u_{12} + u_{30} + 2\cdot u_{00} & > & 6 \cdot u_{11} &&
      3\cdot u_{21} + u_{03} + 2\cdot u_{00} & > & 6 \cdot u_{11}\\
      2\cdot u_{30} + 3\cdot u_{02} + u_{00} & > & 6 \cdot u_{11} &&
      2\cdot u_{03} + 3\cdot u_{20} + u_{00} & > & 6 \cdot u_{11}\\
      u_{12} + u_{30} + u_{00} + u_{02} & > & 4 \cdot u_{11} &&
      u_{21} + u_{03} + u_{00} + u_{20} & > & 4 \cdot u_{11}\\
      u_{01} + u_{10} + u_{03} + u_{30} & > & 4 \cdot u_{11} &&
      2\cdot u_{01} + u_{12} + u_{30} & > & 4 \cdot u_{11}\\
      2\cdot u_{10} + u_{21} + u_{03} & > & 4 \cdot u_{11} &&
      2\cdot u_{12} + u_{20} + u_{00} & > & 4 \cdot u_{11}\\
      2\cdot u_{21} + u_{02} + u_{00} & > & 4 \cdot u_{11} &&
      2\cdot u_{02} + u_{10} + u_{30} & > & 4 \cdot u_{11}\\
      2\cdot u_{20} + u_{01} + u_{03} & > & 4 \cdot u_{11} &&
      u_{20} + u_{01} + u_{12} & > & 3 \cdot u_{11}\\
      u_{02} + u_{10} + u_{21} & > & 3 \cdot u_{11} &&
      u_{30} + u_{02} + u_{01} & > & 3 \cdot u_{11}\\
      u_{03} + u_{10} + u_{20} & > & 3 \cdot u_{11} &&
      u_{00} + u_{12} + u_{21} & > & 3 \cdot u_{11}\\
      u_{00} + u_{30} + u_{03} & > & 3 \cdot u_{11} &&
      u_{21} + u_{01} & > & 2 \cdot u_{11}\\
      u_{10} + u_{12} & > & 2 \cdot u_{11} &&
      u_{20} + u_{02} & > & 2 \cdot u_{11}      
    \end{array}
  \end{equation}
  These inequalities determine thus a cone in $\R^{\ca_3}$ which
  contains $U$. A simple computation with \texttt{polymake}
  allows to compute the extreme rays of this cone and to check that
  they actually satisfy the inequalities of the single cone of the
  Gr\"obner fan of $A$ which is dual to the vertex $12 e_{11}$ in $N(A)$. 
  This proves the claim.
\end{proof}

\begin{remark}
  The computations were done with the procedure \texttt{testInteriorInequalities}
  in the library \texttt{jinvariant.lib} (see \cite{KMM07a}).
  
  The inequalities in \eqref{eq-inequalities} actually determine
  precisely the cone $U$, which is less obvious than that they are
  necessary, but this can again be easily tested using
  \texttt{polymake}. 
\end{remark}
                                %

\section{Rays of the Secondary Fan}\label{sec-rays}

In this section, we classify the rays of the secondary fan of $\ca_3$.
This is part of our second proof of the main theorem.  

\begin{definition} 
  Let $\nu\in\ca_3$ be a lattice point that is
  not a vertex.  The {\em lift} associated to $\nu$ is the ray in
  $\R^{\ca_3}$ consisting of all functions $\psi:{\ca_3}\longrightarrow\R$ of the form
  \[\psi(\omega)=
  \begin{cases}
    a+v\cdot\omega+b, & \text{if }\omega=\nu, \\
    a+v\cdot\omega, & \text{else},
  \end{cases}\]
  for some $a\in\R$, $b\in\R^+$, and $v\in\R^2$.  The marked subdivision
  of $({Q_3},{\ca_3})$ associated to the lift is $\{({Q_3},{\ca_3}\setminus\{\nu\})\}$.
\end{definition}

\begin{figure}[h]
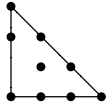

  \begin{texdraw}
    \drawdim cm \relunitscale 0.4 \linewd 0.05 \lpatt (1 0) 
    \move (0 0) \lvec (3 0) \lvec (0 3) \lvec (0 0) 
    \move (0 0) \fcir f:0 r:0.15 
    \move (0 2) \fcir f:0 r:0.15 
    \move (0 3) \fcir f:0 r:0.15 
    \move (1 0) \fcir f:0 r:0.15 
    \move (1 1) \fcir f:0 r:0.15 
    \move (1 2) \fcir f:0 r:0.15
    \move (2 0) \fcir f:0 r:0.15 
    \move (2 1) \fcir f:0 r:0.15 
    \move (3 0) \fcir f:0 r:0.15 
  \end{texdraw}\centering
  \caption{The subdivision associated to the lift with $\nu=(0,1)$.}
\end{figure}

\begin{definition} 
  Let $w\in\R^2$ and $c\in\R$ be such that
  $\{\omega\in\R^2\;|\;\omega\cdot w=c\}$ is a line $l$  through ${Q_3}$
  that intersects the boundary of ${Q_3}$ in lattice points.  The {\em fold through $l$} is the
  cone consisting 
  of all functions $\psi:{\ca_3}\longrightarrow\R$ of the form
  \[\psi(\omega)=a+v\cdot\omega+b\cdot\max\{\omega\cdot w-c,0\}\]
  for some $a\in\R$, $b\in\R^+$, and $v\in\R^2$.
  The subdivision associated to the fold is
  $\{(Q_+,\ca_+),(Q_-,\ca_-)\}$ where
  \begin{eqnarray*}
    Q_+&=&\{v\in {Q_3}|v\cdot w\geq c\}\\
    Q_-&=&\{v\in {Q_3}|v\cdot w\leq c\}
  \end{eqnarray*}
  and $\ca_+=Q_+\cap{\ca_3}$, $\ca_-=Q_-\cap{\ca_3}$.
\end{definition}

\begin{figure}[h]
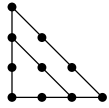

  \begin{texdraw}
  \drawdim cm  \relunitscale 0.4
       \linewd 0.05
        \move (3 0)
        \lvec (0 0)
        \move (0 0)
        \lvec (0 3)
        \move (0 3)
        \lvec (3 0)

        \move (2 0)
        \lvec (0 2)
        \move (0 0) \fcir f:0 r:0.15
        \move (0 1) \fcir f:0 r:0.15
        \move (0 2) \fcir f:0 r:0.15
        \move (0 3) \fcir f:0 r:0.15
        \move (1 0) \fcir f:0 r:0.15
        \move (1 1) \fcir f:0 r:0.15
        \move (1 2) \fcir f:0 r:0.15

        \move (2 0) \fcir f:0 r:0.15
        \move (2 1) \fcir f:0 r:0.15
        \move (3 0) \fcir f:0 r:0.15
  \end{texdraw}\centering
  \caption{A subdivision associated to a fold.}
\end{figure}

\begin{definition} 
  Consider three lattice points $p_1$, $p_2$ and $p_3$ on the boundary of
  $Q_3$ (ordered counterclockwise) and the corresponding vectors
  $v_i$ pointing from the interior point $p=(1,1)$ to $p_i$. The 
  half lines starting in $p$ and passing through $p_i$ divide $\R^2$
  into three cones, say $C_{12}$, $C_{23}$ and $C_{31}$ where
  $C_{ij}$ contains the points $p_i$ and $p_j$ in its boundary. 
  Note that the angle between $v_i$ and $v_{i+1}$ is less than $180$
  degrees.
  Define the function $\phi:{\ca_3}\longrightarrow\R$ by
  \[\phi_b(\omega)=
  \begin{cases}0 &\mbox{ if }\omega\in C_{12}\\
    tb & \mbox{ if }\omega=p+sv_2+tv_3\in C_{23}\\
    tb & \mbox{ if }\omega=p+tv_3+sv_1\in C_{31}
  \end{cases}.\]
  for $b\in\R^+$.  In other words $\phi_b$ is $0$ on $C_1$ and is linear on $p+v_3\cdot\R^+$.
  The {\em pinwheel} through $p_1$, $p_2$ and $p_3$ is the cone
  consisting of all functions $\psi:{\ca_3}\longrightarrow\R$ of the form
  \[\psi(\omega)=a+\omega\cdot v+\phi_b(\omega),\]
  for some $a\in\R$, $b\in\R^+$ and $v\in\R^2$.
  The subdivision associated to the pinwheel is
  subdivision $\{(Q_{12},\ca_{12}),(Q_{23},\ca_{23}),(Q_{31},\ca_{31})\}$
  with $Q_{ij}=C_{ij}\cap Q_3$ and $\ca_{ij}=Q_{ij}\cap\ca_3$.
\end{definition}

\begin{figure}[h]
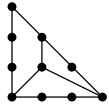

  \begin{texdraw}
       \drawdim cm  \relunitscale 0.4
       \linewd 0.05
        \move (3 0)
        \lvec (0 0)
        \move (0 0)
        \lvec (0 3)
        \move (0 3)
        \lvec (3 0)

        \move (3 0)
        \lvec (1 1)
        \move (1 1)
        \lvec (1 2)
        \move (0 0)
        \lvec (1 1)
        \move (0 0) \fcir f:0 r:0.15
        \move (0 1) \fcir f:0 r:0.15
        \move (0 2) \fcir f:0 r:0.15
        \move (0 3) \fcir f:0 r:0.15
        \move (1 0) \fcir f:0 r:0.15
        \move (1 1) \fcir f:0 r:0.15
        \move (1 2) \fcir f:0 r:0.15

        \move (2 0) \fcir f:0 r:0.15
        \move (2 1) \fcir f:0 r:0.15
        \move (3 0) \fcir f:0 r:0.15
   \end{texdraw}\centering
  \caption{A subdivision associated to a pinwheel.}
\end{figure}

\begin{proposition} \label{prop-rays} 
  Any ray of the secondary fan of ${\ca_3}$ is a lift, a
  fold, or a pinwheel.
\end{proposition}

\begin{proof}
  Note that lifts, folds and pinwheels are obviously rays of the
  secondary fan of $\ca_3$.

  Let $C$ be a ray in the secondary fan of ${\ca_3}$.  Consider the associated
  subdivision $S=\{(Q_\theta,\ca_\theta)\;|\;\theta\in\Theta\}$.  The only
  non-trivial coarsening  of $S$ is the coarsest subdivision.

  If $\bigcup_{\theta\in\Theta} \ca_\theta\neq{\ca_3}$, let
  $\omega\in{\ca_3}\setminus\bigcup_{\theta\in\Theta} \ca_\theta$. 
  The marked subdivision associated to the lift of $\omega$ is a coarsening
  of $S$ and the lift of $\omega$ is a ray of the secondary fan of
  ${\ca_3}$. Therefore the marked subdivision $S$ itself must determine a
  lift. 

  Suppose $\bigcup_{\theta\in\Theta} \ca_\theta={\ca_3}$.  There must be at least $2$ polygons
  in the subdivision.  There must therefore be an edge $E$ of a
  $Q_\theta$ which is not an edge of ${Q_3}$.  If the edge has both
  end-points on the boundary of ${Q_3}$ then the fold through $E$ is a
  coarsening of $S$.  Therefore $S$
  must determine a fold.  Otherwise, one 
  end-point of $E$ must be the interior lattice point, $p$.  There can
  be at most $3$ edges of the subdivision meeting at 
  $p$.  If there were more, we could remove all but $3$ of them
  and still have less than $180$ degrees counterclockwise between
  any two edges.  Therefore, we have a pinwheel.
\end{proof}

There is a natural action of $S_3\subset \PGl_3(\K)$ on
$\Sym^3(\K^3)$.  This induces an action of $S_3$ on the secondary
fan of $\ca_3$.  We list the folds and pinwheels that are rays of the
secondary fan of $\ca_3$ up to $S_3$ action.  

\begin{center}
    \begin{texdraw}
       \drawdim cm  \relunitscale 0.4
       \linewd 0.05
        \move (3 0)
        \lvec (0 0)
        \move (0 0)
        \lvec (0 3)
        \move (0 3)
        \lvec (3 0)

        \move (3 0)
        \lvec (1 1)
        \move (1 1)
        \lvec (1 2)
        \move (0 0)
        \lvec (1 1)
        \move (0 0) \fcir f:0 r:0.15
        \move (0 1) \fcir f:0 r:0.15
        \move (0 2) \fcir f:0 r:0.15
        \move (0 3) \fcir f:0 r:0.15
        \move (1 0) \fcir f:0 r:0.15
        \move (1 1) \fcir f:0 r:0.15
        \move (1 2) \fcir f:0 r:0.15

        \move (2 0) \fcir f:0 r:0.15
        \move (2 1) \fcir f:0 r:0.15
        \move (3 0) \fcir f:0 r:0.15
   \end{texdraw}
   \hspace{1cm}
   \begin{texdraw}
       \drawdim cm  \relunitscale 0.4
       \linewd 0.05
        \move (0 3)
        \lvec (3 0)
        \move (3 0)
        \lvec (0 0)
        \move (0 0)
        \lvec (0 3)

        \move (3 0)
        \lvec (1 1)
        \move (1 1)
        \lvec (0 3)
        \move (0 0)
        \lvec (1 1)
        \move (0 0) \fcir f:0 r:0.15
        \move (0 1) \fcir f:0 r:0.15
        \move (0 2) \fcir f:0 r:0.15
        \move (0 3) \fcir f:0 r:0.15
        \move (1 0) \fcir f:0 r:0.15
        \move (1 1) \fcir f:0 r:0.15
        \move (1 2) \fcir f:0 r:0.15

        \move (2 0) \fcir f:0 r:0.15
        \move (2 1) \fcir f:0 r:0.15
        \move (3 0) \fcir f:0 r:0.15
   \end{texdraw}
   \hspace{1cm}
   \begin{texdraw}
       \drawdim cm  \relunitscale 0.4
       \linewd 0.05
        \move (0 3)
        \lvec (3 0)
        \move (3 0)
        \lvec (0 0)
        \move (0 0)
        \lvec (0 3)

        \move (3 0)
        \lvec (1 1)
        \move (1 1)
        \lvec (0 2)
        \move (0 0)
        \lvec (1 1)
        \move (0 0) \fcir f:0 r:0.15
        \move (0 1) \fcir f:0 r:0.15
        \move (0 2) \fcir f:0 r:0.15
        \move (0 3) \fcir f:0 r:0.15
        \move (1 0) \fcir f:0 r:0.15
        \move (1 1) \fcir f:0 r:0.15
        \move (1 2) \fcir f:0 r:0.15

        \move (2 0) \fcir f:0 r:0.15
        \move (2 1) \fcir f:0 r:0.15
        \move (3 0) \fcir f:0 r:0.15
   \end{texdraw}
   \hspace{1cm}
   \begin{texdraw}
       \drawdim cm  \relunitscale 0.4
       \linewd 0.05
        \move (3 0)
        \lvec (0 0)
        \move (0 0)
        \lvec (0 3)
        \move (0 3)
        \lvec (3 0)

        \move (1 1)
        \lvec (2 1)
        \move (0 0)
        \lvec (1 1)
        \move (1 1)
        \lvec (1 2)
        \move (0 0) \fcir f:0 r:0.15
        \move (0 1) \fcir f:0 r:0.15
        \move (0 2) \fcir f:0 r:0.15
        \move (0 3) \fcir f:0 r:0.15
        \move (1 0) \fcir f:0 r:0.15
        \move (1 1) \fcir f:0 r:0.15
        \move (1 2) \fcir f:0 r:0.15

        \move (2 0) \fcir f:0 r:0.15
        \move (2 1) \fcir f:0 r:0.15
        \move (3 0) \fcir f:0 r:0.15
   \end{texdraw}
   \hspace{1cm}
   \begin{texdraw}
       \drawdim cm  \relunitscale 0.4
       \linewd 0.05
        \move (3 0)
        \lvec (0 0)
        \move (0 0)
        \lvec (0 3)
        \move (0 3)
        \lvec (3 0)

        \move (1 1)
        \lvec (2 1)
        \move (0 0)
        \lvec (1 1)
        \move (1 1)
        \lvec (0 2)
        \move (0 0) \fcir f:0 r:0.15
        \move (0 1) \fcir f:0 r:0.15
        \move (0 2) \fcir f:0 r:0.15
        \move (0 3) \fcir f:0 r:0.15
        \move (1 0) \fcir f:0 r:0.15
        \move (1 1) \fcir f:0 r:0.15
        \move (1 2) \fcir f:0 r:0.15

        \move (2 0) \fcir f:0 r:0.15
        \move (2 1) \fcir f:0 r:0.15
        \move (3 0) \fcir f:0 r:0.15
   \end{texdraw}
 \end{center}
 \begin{center}
   \begin{texdraw}
       \drawdim cm  \relunitscale 0.4
       \linewd 0.05
        \move (3 0)
        \lvec (0 0)
        \move (0 0)
        \lvec (0 3)
        \move (0 3)
        \lvec (3 0)

        \move (1 1)
        \lvec (2 1)
        \move (1 0)
        \lvec (1 1)
        \move (1 1)
        \lvec (0 2)
        \move (0 0) \fcir f:0 r:0.15
        \move (0 1) \fcir f:0 r:0.15
        \move (0 2) \fcir f:0 r:0.15
        \move (0 3) \fcir f:0 r:0.15
        \move (1 0) \fcir f:0 r:0.15
        \move (1 1) \fcir f:0 r:0.15
        \move (1 2) \fcir f:0 r:0.15

        \move (2 0) \fcir f:0 r:0.15
        \move (2 1) \fcir f:0 r:0.15
        \move (3 0) \fcir f:0 r:0.15
   \end{texdraw}
   \hspace{1cm}
   \begin{texdraw}
       \drawdim cm  \relunitscale 0.4
       \linewd 0.05
        \move (0 3)
        \lvec (3 0)
        \move (3 0)
        \lvec (0 0)
        \move (0 0)
        \lvec (0 3)

        \move (3 0)
        \lvec (0 1)
        \move (0 0) \fcir f:0 r:0.15
        \move (0 1) \fcir f:0 r:0.15
        \move (0 2) \fcir f:0 r:0.15
        \move (0 3) \fcir f:0 r:0.15
        \move (1 0) \fcir f:0 r:0.15
        \move (1 1) \fcir f:0 r:0.15
        \move (1 2) \fcir f:0 r:0.15

        \move (2 0) \fcir f:0 r:0.15
        \move (2 1) \fcir f:0 r:0.15
        \move (3 0) \fcir f:0 r:0.15
   \end{texdraw}
   \hspace{1cm}
   \begin{texdraw}
       \drawdim cm  \relunitscale 0.4
       \linewd 0.05
        \move (3 0)
        \lvec (0 0)
        \move (0 0)
        \lvec (0 3)
        \move (0 3)
        \lvec (3 0)

        \move (2 0)
        \lvec (0 1)
        \move (0 0) \fcir f:0 r:0.15
        \move (0 1) \fcir f:0 r:0.15
        \move (0 2) \fcir f:0 r:0.15
        \move (0 3) \fcir f:0 r:0.15
        \move (1 0) \fcir f:0 r:0.15
        \move (1 1) \fcir f:0 r:0.15
        \move (1 2) \fcir f:0 r:0.15

        \move (2 0) \fcir f:0 r:0.15
        \move (2 1) \fcir f:0 r:0.15
        \move (3 0) \fcir f:0 r:0.15
   \end{texdraw}
   \hspace{1cm}
   \begin{texdraw}
       \drawdim cm  \relunitscale 0.4
       \linewd 0.05
        \move (3 0)
        \lvec (0 0)
        \move (0 0)
        \lvec (0 3)
        \move (0 3)
        \lvec (3 0)

        \move (1 0)
        \lvec (0 1)
        \move (0 0) \fcir f:0 r:0.15
        \move (0 1) \fcir f:0 r:0.15
        \move (0 2) \fcir f:0 r:0.15
        \move (0 3) \fcir f:0 r:0.15
        \move (1 0) \fcir f:0 r:0.15
        \move (1 1) \fcir f:0 r:0.15
        \move (1 2) \fcir f:0 r:0.15

        \move (2 0) \fcir f:0 r:0.15
        \move (2 1) \fcir f:0 r:0.15
        \move (3 0) \fcir f:0 r:0.15
   \end{texdraw}
   \hspace{1cm}
   \begin{texdraw}
       \drawdim cm  \relunitscale 0.4
       \linewd 0.05
        \move (3 0)
        \lvec (0 0)
        \move (0 0)
        \lvec (0 3)
        \move (0 3)
        \lvec (3 0)

        \move (2 0)
        \lvec (0 2)
        \move (0 0) \fcir f:0 r:0.15
        \move (0 1) \fcir f:0 r:0.15
        \move (0 2) \fcir f:0 r:0.15
        \move (0 3) \fcir f:0 r:0.15
        \move (1 0) \fcir f:0 r:0.15
        \move (1 1) \fcir f:0 r:0.15
        \move (1 2) \fcir f:0 r:0.15

        \move (2 0) \fcir f:0 r:0.15
        \move (2 1) \fcir f:0 r:0.15
        \move (3 0) \fcir f:0 r:0.15
   \end{texdraw}
 \end{center}


\end{document}

%% file: Graphics/exmarksub.pstex_t
\begin{picture}(0,0)%
\includegraphics{Graphics/exmarksub.pstex}%
\end{picture}%
\setlength{\unitlength}{4144sp}%
\begingroup\makeatletter\ifx\SetFigFont\undefined%
\gdef\SetFigFont#1#2#3#4#5{%
  \reset@font\fontsize{#1}{#2pt}%
  \fontfamily{#3}\fontseries{#4}\fontshape{#5}%
  \selectfont}%
\fi\endgroup%
\begin{picture}(582,582)(5380,-3682)
\end{picture}%

%% file: Graphics/cyclex.pstex_t
\begin{picture}(0,0)%
\includegraphics{Graphics/cyclex.pstex}%
\end{picture}%
\setlength{\unitlength}{4144sp}%
\begingroup\makeatletter\ifx\SetFigFont\undefined%
\gdef\SetFigFont#1#2#3#4#5{%
  \reset@font\fontsize{#1}{#2pt}%
  \fontfamily{#3}\fontseries{#4}\fontshape{#5}%
  \selectfont}%
\fi\endgroup%
\begin{picture}(2499,1520)(4939,-4618)
\end{picture}%

%% file: Graphics/cyclelength.pstex_t
\begin{picture}(0,0)%
\includegraphics{Graphics/cyclelength.pstex}%
\end{picture}%
\setlength{\unitlength}{4144sp}%
\begingroup\makeatletter\ifx\SetFigFont\undefined%
\gdef\SetFigFont#1#2#3#4#5{%
  \reset@font\fontsize{#1}{#2pt}%
  \fontfamily{#3}\fontseries{#4}\fontshape{#5}%
  \selectfont}%
\fi\endgroup%
\begin{picture}(2513,1059)(5105,-3763)
\end{picture}%

%% file: Graphics/cyclelimit.pstex_t
\begin{picture}(0,0)%
\includegraphics{Graphics/cyclelimit.pstex}%
\end{picture}%
\setlength{\unitlength}{4144sp}%
\begingroup\makeatletter\ifx\SetFigFont\undefined%
\gdef\SetFigFont#1#2#3#4#5{%
  \reset@font\fontsize{#1}{#2pt}%
  \fontfamily{#3}\fontseries{#4}\fontshape{#5}%
  \selectfont}%
\fi\endgroup%
\begin{picture}(2769,1295)(4804,-4573)
\put(5086,-4426){\makebox(0,0)[lb]{\smash{{\SetFigFont{10}{12.0}{\familydefault}{\mddefault}{\updefault}{\color[rgb]{0,0,0}$E$}%
}}}}
\put(6706,-4426){\makebox(0,0)[lb]{\smash{{\SetFigFont{10}{12.0}{\familydefault}{\mddefault}{\updefault}{\color[rgb]{0,0,0}$E_2$}%
}}}}
\put(6706,-4066){\makebox(0,0)[lb]{\smash{{\SetFigFont{10}{12.0}{\familydefault}{\mddefault}{\updefault}{\color[rgb]{0,0,0}$E_1$}%
}}}}
\end{picture}%

%% file: Graphics/extrans.pstex_t
\begin{picture}(0,0)%
\includegraphics{Graphics/extrans.pstex}%
\end{picture}%
\setlength{\unitlength}{4144sp}%
\begingroup\makeatletter\ifx\SetFigFont\undefined%
\gdef\SetFigFont#1#2#3#4#5{%
  \reset@font\fontsize{#1}{#2pt}%
  \fontfamily{#3}\fontseries{#4}\fontshape{#5}%
  \selectfont}%
\fi\endgroup%
\begin{picture}(4080,1966)(2371,-4670)
\put(5176,-2851){\makebox(0,0)[lb]{\smash{{\SetFigFont{10}{12.0}{\familydefault}{\mddefault}{\updefault}{\color[rgb]{0,0,0}$u_{03}$}%
}}}}
\put(5806,-3256){\makebox(0,0)[lb]{\smash{{\SetFigFont{10}{12.0}{\familydefault}{\mddefault}{\updefault}{\color[rgb]{0,0,0}$u_{12}$}%
}}}}
\put(6301,-3481){\makebox(0,0)[lb]{\smash{{\SetFigFont{10}{12.0}{\familydefault}{\mddefault}{\updefault}{\color[rgb]{0,0,0}$\min\{u_{11},u_{21}+b\}$}%
}}}}
\put(6256,-3796){\makebox(0,0)[lb]{\smash{{\SetFigFont{10}{12.0}{\familydefault}{\mddefault}{\updefault}{\color[rgb]{0,0,0}$u_{21}$}%
}}}}
\put(6436,-4066){\makebox(0,0)[lb]{\smash{{\SetFigFont{10}{12.0}{\familydefault}{\mddefault}{\updefault}{\color[rgb]{0,0,0}$u_{30}$}%
}}}}
\put(5941,-4381){\makebox(0,0)[lb]{\smash{{\SetFigFont{10}{12.0}{\familydefault}{\mddefault}{\updefault}{\color[rgb]{0,0,0}$\min\{u_{20},u_{30}+b\}$}%
}}}}
\put(5266,-4606){\makebox(0,0)[lb]{\smash{{\SetFigFont{10}{12.0}{\familydefault}{\mddefault}{\updefault}{\color[rgb]{0,0,0}$\min\{u_{10},u_{20}+b,u_{30}+2b\}$}%
}}}}
\put(2386,-4066){\makebox(0,0)[lb]{\smash{{\SetFigFont{10}{12.0}{\familydefault}{\mddefault}{\updefault}{\color[rgb]{0,0,0}$\min\{u_{00},u_{10}+b,u_{20}+2b,u_{30}+3b\}$}%
}}}}
\put(3016,-3796){\makebox(0,0)[lb]{\smash{{\SetFigFont{10}{12.0}{\familydefault}{\mddefault}{\updefault}{\color[rgb]{0,0,0}$\min\{u_{01},u_{11}+b,u_{21}+2b\}$}%
}}}}
\put(3646,-3526){\makebox(0,0)[lb]{\smash{{\SetFigFont{10}{12.0}{\familydefault}{\mddefault}{\updefault}{\color[rgb]{0,0,0}$\min\{u_{02},u_{12}+b\}$}%
}}}}
\end{picture}%

%% file: Graphics/raysec.pstex_t
\begin{picture}(0,0)%
\includegraphics{Graphics/raysec.pstex}%
\end{picture}%
\setlength{\unitlength}{4144sp}%
\begingroup\makeatletter\ifx\SetFigFont\undefined%
\gdef\SetFigFont#1#2#3#4#5{%
  \reset@font\fontsize{#1}{#2pt}%
  \fontfamily{#3}\fontseries{#4}\fontshape{#5}%
  \selectfont}%
\fi\endgroup%
\begin{picture}(594,594)(4474,-3508)
\end{picture}%

%% file: Graphics/deltaeq.pstex_t
\begin{picture}(0,0)%
\includegraphics{Graphics/deltaeq.pstex}%
\end{picture}%
\setlength{\unitlength}{4144sp}%
\begingroup\makeatletter\ifx\SetFigFont\undefined%
\gdef\SetFigFont#1#2#3#4#5{%
  \reset@font\fontsize{#1}{#2pt}%
  \fontfamily{#3}\fontseries{#4}\fontshape{#5}%
  \selectfont}%
\fi\endgroup%
\begin{picture}(1554,564)(5209,-3763)
\end{picture}%

%% file: Graphics/circuit.pstex_t
\begin{picture}(0,0)%
\includegraphics{Graphics/circuit.pstex}%
\end{picture}%
\setlength{\unitlength}{4144sp}%
\begingroup\makeatletter\ifx\SetFigFont\undefined%
\gdef\SetFigFont#1#2#3#4#5{%
  \reset@font\fontsize{#1}{#2pt}%
  \fontfamily{#3}\fontseries{#4}\fontshape{#5}%
  \selectfont}%
\fi\endgroup%
\begin{picture}(3225,1381)(5476,-5165)
\put(6121,-5101){\makebox(0,0)[lb]{\smash{{\SetFigFont{12}{14.4}{\familydefault}{\mddefault}{\updefault}{\color[rgb]{0,0,0}$T$}%
}}}}
\put(6121,-4201){\makebox(0,0)[lb]{\smash{{\SetFigFont{10}{12.0}{\familydefault}{\mddefault}{\updefault}{\color[rgb]{0,0,0}$\omega_i$}%
}}}}
\put(5716,-3931){\makebox(0,0)[lb]{\smash{{\SetFigFont{10}{12.0}{\familydefault}{\mddefault}{\updefault}{\color[rgb]{0,0,0}$\omega_{i-1}$}%
}}}}
\put(6436,-4561){\makebox(0,0)[lb]{\smash{{\SetFigFont{10}{12.0}{\familydefault}{\mddefault}{\updefault}{\color[rgb]{0,0,0}$\omega_{i+1}$}%
}}}}
\put(5491,-4246){\makebox(0,0)[lb]{\smash{{\SetFigFont{10}{12.0}{\familydefault}{\mddefault}{\updefault}{\color[rgb]{0,0,0}$\omega_{i-2}$}%
}}}}
\put(6121,-4741){\makebox(0,0)[lb]{\smash{{\SetFigFont{10}{12.0}{\familydefault}{\mddefault}{\updefault}{\color[rgb]{0,0,0}$\omega_{i+2}$}%
}}}}
\put(7966,-3931){\makebox(0,0)[lb]{\smash{{\SetFigFont{10}{12.0}{\familydefault}{\mddefault}{\updefault}{\color[rgb]{0,0,0}$\omega_{i-1}$}%
}}}}
\put(8686,-4561){\makebox(0,0)[lb]{\smash{{\SetFigFont{10}{12.0}{\familydefault}{\mddefault}{\updefault}{\color[rgb]{0,0,0}$\omega_{i+1}$}%
}}}}
\put(7741,-4246){\makebox(0,0)[lb]{\smash{{\SetFigFont{10}{12.0}{\familydefault}{\mddefault}{\updefault}{\color[rgb]{0,0,0}$\omega_{i-2}$}%
}}}}
\put(8371,-4741){\makebox(0,0)[lb]{\smash{{\SetFigFont{10}{12.0}{\familydefault}{\mddefault}{\updefault}{\color[rgb]{0,0,0}$\omega_{i+2}$}%
}}}}
\put(8416,-5101){\makebox(0,0)[lb]{\smash{{\SetFigFont{12}{14.4}{\familydefault}{\mddefault}{\updefault}{\color[rgb]{0,0,0}$T'$}%
}}}}
\put(5896,-4516){\makebox(0,0)[lb]{\smash{{\SetFigFont{10}{12.0}{\familydefault}{\mddefault}{\updefault}{\color[rgb]{0,0,0}$\omega$}%
}}}}
\put(8146,-4516){\makebox(0,0)[lb]{\smash{{\SetFigFont{10}{12.0}{\familydefault}{\mddefault}{\updefault}{\color[rgb]{0,0,0}$\omega$}%
}}}}
\end{picture}%

%% file: Graphics/L1.pstex_t
\begin{picture}(0,0)%
\includegraphics{Graphics/L1.pstex}%
\end{picture}%
\setlength{\unitlength}{4144sp}%
\begingroup\makeatletter\ifx\SetFigFont\undefined%
\gdef\SetFigFont#1#2#3#4#5{%
  \reset@font\fontsize{#1}{#2pt}%
  \fontfamily{#3}\fontseries{#4}\fontshape{#5}%
  \selectfont}%
\fi\endgroup%
\begin{picture}(1658,578)(4752,-3770)
\end{picture}%

%% file: Graphics/L2.pstex_t
\begin{picture}(0,0)%
\includegraphics{Graphics/L2.pstex}%
\end{picture}%
\setlength{\unitlength}{4144sp}%
\begingroup\makeatletter\ifx\SetFigFont\undefined%
\gdef\SetFigFont#1#2#3#4#5{%
  \reset@font\fontsize{#1}{#2pt}%
  \fontfamily{#3}\fontseries{#4}\fontshape{#5}%
  \selectfont}%
\fi\endgroup%
\begin{picture}(1558,566)(4757,-3764)
\end{picture}%

%% file: Graphics/L3.pstex_t
\begin{picture}(0,0)%
\includegraphics{Graphics/L3.pstex}%
\end{picture}%
\setlength{\unitlength}{4144sp}%
\begingroup\makeatletter\ifx\SetFigFont\undefined%
\gdef\SetFigFont#1#2#3#4#5{%
  \reset@font\fontsize{#1}{#2pt}%
  \fontfamily{#3}\fontseries{#4}\fontshape{#5}%
  \selectfont}%
\fi\endgroup%
\begin{picture}(1564,574)(4844,-3588)
\end{picture}%

%% file: Graphics/L4.pstex_t
\begin{picture}(0,0)%
\includegraphics{Graphics/L4.pstex}%
\end{picture}%
\setlength{\unitlength}{4144sp}%
\begingroup\makeatletter\ifx\SetFigFont\undefined%
\gdef\SetFigFont#1#2#3#4#5{%
  \reset@font\fontsize{#1}{#2pt}%
  \fontfamily{#3}\fontseries{#4}\fontshape{#5}%
  \selectfont}%
\fi\endgroup%
\begin{picture}(3276,574)(4483,-3588)
\end{picture}%

%% file: Graphics/L5.pstex_t
\begin{picture}(0,0)%
\includegraphics{Graphics/L5.pstex}%
\end{picture}%
\setlength{\unitlength}{4144sp}%
\begingroup\makeatletter\ifx\SetFigFont\undefined%
\gdef\SetFigFont#1#2#3#4#5{%
  \reset@font\fontsize{#1}{#2pt}%
  \fontfamily{#3}\fontseries{#4}\fontshape{#5}%
  \selectfont}%
\fi\endgroup%
\begin{picture}(1484,582)(5199,-3952)
\end{picture}%

%% file: Graphics/L6.pstex_t
\begin{picture}(0,0)%
\includegraphics{Graphics/L6.pstex}%
\end{picture}%
\setlength{\unitlength}{4144sp}%
\begingroup\makeatletter\ifx\SetFigFont\undefined%
\gdef\SetFigFont#1#2#3#4#5{%
  \reset@font\fontsize{#1}{#2pt}%
  \fontfamily{#3}\fontseries{#4}\fontshape{#5}%
  \selectfont}%
\fi\endgroup%
\begin{picture}(1486,586)(5198,-3774)
\end{picture}%

%% file: Graphics/L7.pstex_t
\begin{picture}(0,0)%
\includegraphics{Graphics/L7.pstex}%
\end{picture}%
\setlength{\unitlength}{4144sp}%
\begingroup\makeatletter\ifx\SetFigFont\undefined%
\gdef\SetFigFont#1#2#3#4#5{%
  \reset@font\fontsize{#1}{#2pt}%
  \fontfamily{#3}\fontseries{#4}\fontshape{#5}%
  \selectfont}%
\fi\endgroup%
\begin{picture}(2388,586)(5107,-3774)
\end{picture}%

%% file: Graphics/L8.pstex_t
\begin{picture}(0,0)%
\includegraphics{Graphics/L8.pstex}%
\end{picture}%
\setlength{\unitlength}{4144sp}%
\begingroup\makeatletter\ifx\SetFigFont\undefined%
\gdef\SetFigFont#1#2#3#4#5{%
  \reset@font\fontsize{#1}{#2pt}%
  \fontfamily{#3}\fontseries{#4}\fontshape{#5}%
  \selectfont}%
\fi\endgroup%
\begin{picture}(583,581)(5199,-3682)
\end{picture}%

%% file: Graphics/boundary.pstex_t
\begin{picture}(0,0)%
\includegraphics{Graphics/boundary.pstex}%
\end{picture}%
\setlength{\unitlength}{4144sp}%
\begingroup\makeatletter\ifx\SetFigFont\undefined%
\gdef\SetFigFont#1#2#3#4#5{%
  \reset@font\fontsize{#1}{#2pt}%
  \fontfamily{#3}\fontseries{#4}\fontshape{#5}%
  \selectfont}%
\fi\endgroup%
\begin{picture}(5094,594)(5104,-3868)
\end{picture}%

%% file: Graphics/vertex.pstex_t
\begin{picture}(0,0)%
\includegraphics{Graphics/vertex.pstex}%
\end{picture}%
\setlength{\unitlength}{4144sp}%
\begingroup\makeatletter\ifx\SetFigFont\undefined%
\gdef\SetFigFont#1#2#3#4#5{%
  \reset@font\fontsize{#1}{#2pt}%
  \fontfamily{#3}\fontseries{#4}\fontshape{#5}%
  \selectfont}%
\fi\endgroup%
\begin{picture}(5090,590)(4926,-3686)
\end{picture}%

%% file: Graphics/Lvert.pstex_t
\begin{picture}(0,0)%
\includegraphics{Graphics/Lvert.pstex}%
\end{picture}%
\setlength{\unitlength}{4144sp}%
\begingroup\makeatletter\ifx\SetFigFont\undefined%
\gdef\SetFigFont#1#2#3#4#5{%
  \reset@font\fontsize{#1}{#2pt}%
  \fontfamily{#3}\fontseries{#4}\fontshape{#5}%
  \selectfont}%
\fi\endgroup%
\begin{picture}(2392,590)(5285,-5846)
\end{picture}%